\documentclass[12pt]{amsart}
\usepackage{fullpage}

\usepackage{hyperref}

\newcommand{\ord}{\mbox{ord}}

\usepackage{color, graphicx}
\usepackage{amsmath}
\usepackage{amssymb}

\usepackage{charter}
\usepackage[T1]{fontenc}
\usepackage{textcomp}

\usepackage{calrsfs}

\theoremstyle{plain}

\newtheorem{theorem}{Theorem}[section]
\newtheorem{lemma}[theorem]{\bf Lemma}
\newtheorem{corollary}[theorem]{\bf Corollary}
\newtheorem{proposition}[theorem]{\bf Proposition}

\theoremstyle{definition}
\newtheorem{definition}[theorem]{\bf Definition}

\theoremstyle{remark}
\newtheorem{remark}[theorem]{\bf Remark}

\newtheorem{notation}[theorem]{\bf Notation}
\newtheorem{notationassumptions}[theorem]{\bf Notation and Assumptions}

\newcommand{\bad}{{\operatorname{bad}}}
\newcommand{\tors}{{\operatorname{tors}}}
\newcommand{\calA}{{\mathcal A}}
\newcommand{\calB}{{\mathcal B}}
\newcommand{\calC}{{\mathcal C}}

\newcommand{\calP}{{\mathcal P}}

\newcommand{\calS}{{\mathcal S}}

\newcommand{\calV}{{\mathcal V}}
\newcommand{\calW}{{\mathcal W}}
\newcommand{\calX}{{\mathcal X}}
\newcommand{\calY}{{\mathcal Y}}

\newcommand{\F}{{\mathbb F}}

\newcommand{\Q}{{\mathbb Q}}

\newcommand{\Z}{{\mathbb Z}}

\newcommand{\pp}{{\mathfrak p}}
\newcommand{\qq}{{\mathfrak q}}

\newcommand{\nn}{{\mathfrak n}}
\newcommand{\dd}{{\mathfrak d}}
\newcommand{\cc}{{\mathfrak c}}
\newcommand{\ttt}{{\mathfrak t}}

  \newcommand{\ff}{{\mathfrak f}}
  \newcommand{\hh}{{\mathfrak h}}
\newcommand{\be}{\begin{enumerate}}
\newcommand{\ee}{\end{enumerate}}
\newcommand{\norm}{{\mathbf N}}

\def\ord{\mathop{\mathrm{ord}}\nolimits}

\usepackage[OT2,OT1]{fontenc}
\newcommand\cyr{%
\renewcommand\rmdefault{wncyr}%
\renewcommand\sfdefault{wncyss}%
\renewcommand\encodingdefault{OT2}%
\normalfont
\selectfont}
\DeclareTextFontCommand{\textcyr}{\cyr}

\begin{document}
\bibliographystyle{plain}%
 \title{Using Indices of Points on an Elliptic Curve to Construct A Diophantine Model of $\Z$ and Define $\Z$ Using One
Universal Quantifier in Very Large Subrings of Number Fields, Including $\Q$ }%
\author{Alexandra Shlapentokh}%
\thanks{The research for this paper has been partially supported by NSF grant DMS-0650927,  a grant from John Templeton Foundation and ECU Faculty Senate Summer 2011 grant. The author would also like to thank Bjorn Poonen for his help.}
\address{Department of Mathematics \\ East Carolina University \\ Greenville, NC 27858}%
\email{shlapentokha@ecu.edu }
\urladdr{www.personal.ecu.edu/shlapentokha} 

\date{\today}

\maketitle

\bigskip

\bigskip

\begin{abstract}
Let $K$ be a number field such that there exists an elliptic curve  $E$ of rank one over $K$. For a set $\calW_K$
 of primes of $K$, let $O_{K,\calW_K}=\{x\in K: \ord_{\pp}x \geq 0, \, \forall \pp \not \in \calW_K\}$. Let $P \in E(K)$
be a generator of $E(K)$ modulo the torsion subgroup.  Let $(x_n(P),y_n(P))$ be the affine coordinates of $[n]P$ with respect to a fixed Weierstrass
equation of $E$. We show that there exists a set $\calW_K$ of primes of $K$ of natural  density one such that in
$O_{K,\calW_K}$ multiplication of indices (with respect to some fixed multiple of $P$) is existentially definable and therefore these indices can be
used to construct a Diophantine model of $\Z$.  We also show that $\Z$ is definable over $O_{K,\calW_K}$ using just one
universal quantifier.  Both, the construction of a Diophantine model using the indices and the first-order definition of $\Z$
can be lifted to the integral closure of $O_{K,\calW_K}$ in any infinite extension $K_{\infty}$ of $K$ as long
as $E(K_{\infty})$ is finitely generated and of rank one.

\end{abstract}%

\section{Introduction}
The interest in constructing Diophantine models of $\Z$ over various rings and related issues of Diophantine decidability and
definability over rings goes back to a question that was posed by Hilbert: given an arbitrary polynomial equation in several
variables over $\Z$, is there a uniform algorithm to determine whether such an equation has solutions in $\Z$? This
question, otherwise known as Hilbert's Tenth Problem, has been answered negatively in the work of
M. Davis, H. Putnam, J. Robinson and Yu. Matijasevich. (See \cite{Da1}, \cite{Da2} and
\cite{Mate}.) Since the time when this result was obtained, similar questions have been raised for
other fields and rings. In other words, if $R$ is a recursive ring, then, given an arbitrary
polynomial equation in several variables over $R$, is there a uniform algorithm to determine
whether such an equation has solutions in $R$? One way to resolve the question of Diophantine
decidability negatively over a ring of characteristic 0 is to construct a Diophantine definition of
$\Z$ over such a ring. This notion is defined below.

\begin{definition}
Let $R$ be a ring and let $A \subset R^k, k \in \Z_{>0}$.  Then we say that $A$ has a Diophantine definition over
$R$ if there exists a polynomial
\[%
f(t_1,\ldots,t_k,x_1,\ldots,x_n) \in R[t_1,\ldots,t_k,x_1,\ldots,x_n]
\]%
 such that for any $\bar t \in R^k$,
\[
\exists x_1,\ldots,x_n \in R, f(t_1,\ldots,t_k,x_1,...,x_n) = 0 \Longleftrightarrow \bar t \in A.
\]

If the quotient field of $R$ is not algebraically closed,   we can allow
a Diophantine definition to consist of several polynomials without changing the nature of the
relation. (See \cite{Da2} for more details.)
\end{definition}

The usefulness  of Diophantine definitions stems from the following easy lemma.
\begin{lemma}
Let $R_1 \subset R_2$ be two recursive rings such that the quotient field of $R_2$ is not
algebraically closed. Assume that  Hilbert's Tenth Problem (abbreviated as ``HTP'' in the future)   is
undecidable over $R_1$, and $R_1$ has a Diophantine definition over $R_2$. Then HTP is
undecidable over $R_2$.%
\end{lemma}

Using norm equations, Diophantine definitions have been obtained for $\Z$ over the rings of algebraic integers of
some number fields. Jan Denef has constructed a Diophantine definition of $\Z$ for the finite degree totally real
extensions of $\Q$. Jan Denef and Leonard Lipshitz extended Denef's results to all the extensions of degree 2 of
the finite degree totally real fields. Thanases Pheidas and the author of this paper have independently
constructed Diophantine definitions of $\Z$ for number fields with exactly one pair of non-real
embeddings. Finally Harold N. Shapiro and the author of this paper showed that the subfields of all the fields
mentioned above ``inherited'' the Diophantine definitions of $\Z$. (These subfields include all the abelian
extensions.) The proofs of the results listed above can be found in \cite{Den1}, \cite{Den2},
\cite{Den3}, \cite{Ph1}, \cite{Sha-Sh}, and \cite{Sh2}.\\ %

The author modified the norm method to obtain Diophantine definitions of $\Z$ for ``large'' subrings of totally real
number fields (not equal to $\Q$) and their extensions of degree 2.  (See \cite{Sh1}, \cite{Sh6}, \cite{Sh3},
\cite{Sh36}.)  Further, again using norm equations, the author also showed that in some totally real infinite algebraic
extensions of $\Q$ and extensions of degree 2 of such fields one can give  a Diophantine definition of $\Z$ over
integral closures of ``small'' and ``large'' rings, though not over the rings of algebraic integers.  (The terms
``large'' and ''small'' rings will be explained below in Definition \ref{def:large}.)\\

Using elliptic curves Bjorn Poonen has shown  the following in \cite{Po}.%
\begin{theorem}%
\label{thm:po}%
Let $M/K$ be a number field extension with an elliptic curve $E$ defined over $K$, of rank one over $K$, such that
the rank of $E$ over $M$ is also one.  Then $O_K$ (the ring of integers of $K$) is Diophantine over $O_M$.\\%
\end{theorem}%

 Cornelissen, Pheidas and Zahidi weakened somewhat assumptions of Poonen's theorem. Instead of requiring a rank 1
curve retaining its rank in the extension, they require existence of a rank 1 elliptic curve over the bigger field
and an abelian variety over the smaller field retaining its positive rank in the extension (see \cite{CPZ}).  Further,
Poonen and the author have independently shown that the conditions of Theorem \ref{thm:po} can be weakened to
remove the assumption that the rank is one and require only that the rank in the extension is positive and the same as the rank over the ground field
(see \cite{Sh33} and \cite{Po3}).  Following Denef in \cite{Den2}, the author also considered the situations where elliptic curves had finite rank in
infinite extensions and showed that when this happens in a totally real field one can existentially define $\Z$ over the ring of integers of this field
and the ring of integers of any extension of degree 2 of such a field (see \cite{Sh37}).

Recently, in \cite{MR}, Mazur and Rubin showed that if Shafarevich-Tate conjecture held over a number field $K$, then for any prime degree cyclic extension $M$ of $K$, there existed an elliptic curve of rank one over $K$, keeping its rank over $M$.  Combined with Theorem \ref{thm:po}, this new result showed that Shafarevich -Tate conjecture implied HTP is undecidable over the rings of integers of any number field.  Similar consequences can be derived for big rings in any number field.

Perhaps the most prominent open question in the subject is the Diophantine status of $\Q$.  As indicated above, one way to
show unsolvability of HTP over $\Q$ would be to construct a Diophantine definition of $\Z$ over $\Q$.  A Diophantine
definition is an example of a Diophantine model.  Given two recursive rings $R_1$ and $R_2$ we say that $R_2$ has a Diophantine
model of $R_1$ if there exists an injective and recursive map $\phi: R_1 \longrightarrow R_2$ sending Diophantine sets to
Diophantine sets.  If $R_1$ has undecidable Diophantine sets, then so does $R_2$.  Therefore, any recursive ring with a
Diophantine model of $\Z$  has undecidable Diophantine sets and thus HTP is unsolvable over this ring.

It is also not hard to show that given an injection $\phi$ of $\Z$ into a  recursive ring $R$, it is enough to
show that the images of the graphs of addition and multiplication are Diophantine over $R$, in order to conclude that $\phi$
is a Diophantine model.  An old plan for constructing a Diophantine model of $\Z$ over $\Q$ involved elliptic curves of rank
one (see \cite{Ph2}).  More specifically let $E$ be an elliptic curve defined and of rank one over $\Q$.  Fix an affine
Weierstrass equation for $E$, as well as a generator $Q$.  Let $r$ be the size of the torsion group and let $P=[r]Q$.  Let
$(x_n(P),y_n(P))$ be the coordinates of $[n]P$ derived from our fixed affine Weierstrass equation.  Now for $n \not=0$ send
$n$ to $y_n$.  It is easy to see that the graph of addition is Diophantine over $\Q$, but it is not clear what happens to the
graph of multiplication.  This plan has another potentially fatal complication: Mazur's conjectures (see \cite{M1}, \cite{M2},
\cite{M4}).  As was shown in \cite{CZ}, if Mazur's conjecture on topology of rational points holds, there is no Diophantine
model of $\Z$ over $\Q$.   It is precisely these difficulties preventing the resolution of the problem over $\Q$ that
motivated the investigation of Diophantine definability and decidability over ``large'' or ''big'' rings.
These rings can be found in any number field and we define them below.
\begin{definition}%
\label{def:large}
 Let $K$ be a number field and let $\calW_K$ be a set of primes of $K$. Define
$O_{K,\calW_K}$ to be the following ring:
\[%
O_{K,\calW_K} := \{x \in K: \ord_{\pp}x \geq 0, \,\, \forall \pp \not \in \calW_K\}.
\]%
If $\calW_K$ is infinite we will call these rings ``big'' or ``large''.  If $\calW_K$ is finite we refer to the corresponding rings as ``small''.  Such rings are also known as the rings of $\calS$-integers.
\end{definition}%

Perhaps the most significant result concerning big rings was obtained by Poonen in \cite{Po2}.  In this paper he showed that
there exists a big ring inside $\Q$ where the set of primes allowed in the denominator is of natural
density one and the ring possesses a Diophantine model of $\Z$.  To carry out his construction,
Poonen modeled integers by approximation.  More specifically in \cite{Po2} he proved the following.
Let $E$ be a curve of rank one over $\Q$ without complex multiplication and with only one connected component. Let $P$ be a generator of
$E(\Q)$. Then for some set $\calW_{\Q}$ of rational primes of natural density one, we have that
$E(O_{K,\calW_{\Q}}) =\{(x_{\ell_i}, y_{\ell_i}), i \in \Z_{>0}\} \cup \{\mbox{ finite set }\}$,
where $(x_n, y_n)$ are the coordinates of $[n]P$ obtained from a fixed affine Weierstrass equation
of $E$.  Further it is also the case that $|y_{\ell_i} - i| < 10^{-i}$ for all positive integers
$i$.  Later in \cite{PS}, this result was lifted to all number fields with rank one elliptic curves
(also including curves with complex multiplication) though construction of the model proceeded along
a different path but still using a subsequence of coordinates ${(x_{\ell_i}, y_{\ell_i})}$.\\

In this paper we resurrect in a manner of speaking the old plan of modeling $\Z$ using the indices of points on an elliptic
curve but only over a big ring.  More precisely we prove the following theorem.
\begin{theorem}%
\label{thm:main}%
 Let $K$ be a number field.  Let $E$ be an elliptic curve defined and of rank one over $K$.  Let $P$ be a
generator of $E(K)$ modulo the torsion subgroup, and fix an affine Weierstrass equation for $E$ of the form $y^2=x^3+ ax +b$, with $a, b \in O_K$, where $O_K$ is the ring of integers of $K$.  Let $(x_n,y_n)$ be the
coordinates of $[n]P$ derived from this Weierstrass equation.  Then there exists a set of $K$-primes $\calW_K$ of
natural density one,  and a positive integer $m_0$ such that the following set $\Pi \subset O_{K,\calW_K}^{12}$ is Diophantine over $O_{K,\calW_K}$.
\[%
(U_1, U_2, U_3, X_1,X_2,X_3, V_1,V_2, V_3, Y_1, Y_2, Y_3)\in \Pi \Leftrightarrow
\]%
\[%
\exists \mbox{ unique } k_1, k_2, k_3 \in \Z_{\not= 0}   \mbox{ such that } \left(\frac{U_i}{V_i}, \frac{X_i}{Y_i}\right )=(x_{m_0k_i}, y_{m_0k_i}), \mbox{ for } i=1,2,3, \mbox{ and } k_3=k_1k_2.
\]%
\end{theorem}%

We can use this result to construct yet another variation of a Diophantine model of $\Z$.
\begin{definition}%
Let $R$ be a countable recursive ring, let $D \subset R^k, k \in \Z_{>0}$, be a Diophantine subset, and let $\approx$
be a (Diophantine) equivalence relation on $D$, i.e.\ assume that the set $\{(\bar x,\bar y): \bar x,\bar y \in D, \bar x
\approx \bar y\}$ is  a Diophantine subset of $R^{2k}$. Let $D=\bigcup_{i \in \Z}D_i$, where $D_i$ is an
equivalence class of $\approx$, and let $\phi: \Z \longrightarrow \{D_i, i \in \Z\}$ be defined by
$\phi(i) = D_i$. Finally assume that the sets
\[%
Plus=\{(\bar x,\bar y,\bar z): \bar x \in D_i, \bar y\in D_j, \bar z \in D_{i+j}\}
\]%
and
\[%
Times =\{(\bar x,\bar y,\bar z): \bar x \in D_i, \bar y\in D_j, \bar z \in D_{ij}\}
\]%
 are Diophantine over $R$.

Then we will say that $R$ has a \emph{class} Diophantine model of $\Z$.
\end{definition}%

It is clear that if $R$ does
have a class Diophantine model of $\Z$ then HTP is not solvable over $R$. Such a model of $\Z$ has been used already to show
Diophantine undecidability of function fields of positive characteristic (see  \cite{Eis},
\cite{K-R1},   \cite{Ph3},   \cite{Sh30},  \cite{Sh13}, \cite{Sh15}).\\

As a corollary of Theorem \ref{thm:main} we immediately obtain the following statement.
\begin{corollary}%
\label{cor:main}
In the notation above, for $n \not = 0$ let $\phi(n)=[ (U_{m_0n},X_{m_0n}, V_{m_0n}, Y_{m_0n})]$, the equivalence class of $(U_{m_0n},X_{m_0n}, V_{m_0n}, Y_{m_0n})$ under the equivalence relation described below,   where $U_{m_0n},X_{m_0n}$, $V_{m_0n}, Y_{m_0n} \in O_{K,\calW_K}, V_{m_0n}Y_{m_0n} \not = 0$, and
$\displaystyle (x_{m_0n},y_{m_0n})=\left(\frac{U_{m_0n}}{V_{m_0n}},\frac{X_{m_0n}}{Y_{m_0n}} \right)$.  Let $\phi(0)=\{(0,0,0,0)\}$.  Then $\phi$  is a class Diophantine model of $\Z$. (Here  if  $V\hat V\hat Y Y\not=0$ we set $(U,X,V,Y) \approx (\hat{ U},\hat X, \hat V, \hat Y)$ if and only if $\displaystyle \frac{\hat{U}}{\hat{V}}=\frac{U}{V}$ and $\displaystyle \frac{{\hat X}}{{\hat Y}}=\frac{X}{Y}$.)
\end{corollary}%
Using Theorem \ref{thm:main} we also prove the following.
\begin{theorem}%
\label{thm:first-order}%
Let $K$ be a number field. Let $E$ be an elliptic curve defined and of rank one over $K$. Then there exists a set $\calW_K$ of
primes of $K$ of natural density one such that $\Z$ is first-order definable over $O_{K,\calW_K}$ using just one universal
quantifier.
\end{theorem}%
This result is an improvement of  the first-order definability results for big rings in \cite{CS}
and \cite{PO4}, where the first-order definition of $\Z$ was given using just one universal
quantifier over big rings contained in $\Q$ in \cite{PO4} and in some number fields in \cite{CS} with the natural density of the inverted primes
arbitrarily close but not equal to one. (We should also note here that the main result of \cite{PO4} is defining $\Z$ over $\Q$ using two universal
quantifiers.) The result of this paper is also a natural complement to the results of \cite{CZ2} where it was shown that a model of $\Z$ can be defined over $\Q$ using just one universal quantifier provided a certain conjecture on elliptic curves is true.  More recently, Jochen Koenigsmann showed in \cite{Koenig2} that $\Z$ can be defined over $\Q$ using just one universal quantifier.

Finally, Theorem \ref{thm:main} allows us to simplify some results concerning infinite
extensions from \cite{Sh37}.  The result of  Theorem \ref{thm:main}  holds for any algebraic
extension of $\Q$ with a rank 1 finitely generated elliptic curve. No additional assumptions are
required.  In the past we needed some way to define integrality at a prime in an infinite extension to
use this kind of elliptic curve technique.\\

We  finish this section with a notation set to be used in the rest of the paper.

\begin{notation}
\begin{itemize}
\item Let $\calP_{\Q}=\{2,3,5,\dots\}$ denote  the set of rational primes.%
\item Let $K$ be a number field.%
\item Let $\calP_K$ be the set of all finite primes of $K$.%
\item Given $x \in K$,  let $\nn(x) = \prod_{\pp}\pp^{\ord_{\pp}x}$, where the product is taken over
all $\pp \in \calP_K$ such that $\ord_{\pp}x>0$. Let $\dd(x) =\nn(x^{-1})$.%
\item Let $\calW_K \subset \calP_K$ (we will make $\calW_K$ more specific in the next section). \item Let $A, B
\in O_{K,\calW_K}$.  Then we will say that $(A,B)_{\calW_K}=1$ if for all $\pp \in \calP_K \setminus \calW_K$
we have that either $\ord_{\pp}A=0$ or $\ord_{\pp}B=0$.%
\item Let $A, B \in O_{K,\calW_K}$. Then we will say that $A \Big{|}_{\calW_K} B$ if for all $\pp \not \in
\calW_K$ we have that $\ord_{\pp}B \geq \ord_{\pp}A$ or in other words $A$ divides $B$ in the ring
$O_{K,\calW_K}$.
\item Let $h_K$ be the class number of $K$.  (See \cite{Januz}, Chapter I, \S4 for the definition of a class number.)
\item Let $\mathfrak A, \mathfrak B$ be two integral divisors of $K$.  Then we will say that $\mathfrak A \Big{|} \mathfrak B$
to mean that for all $\pp \in \calP(K)$ we have that $\ord_{\pp}\mathfrak A \leq \ord_{\pp}\mathfrak B$. %
\item Suppose $\mathfrak A, \mathfrak B$ are two divisors of $K$ with $\mathfrak B =\mathfrak A^j$.  Then we
set $\sqrt[j]{\mathfrak B}=\mathfrak A$.
\end{itemize}%
 \end{notation}%

\section{An Outline of the Proof of Theorem \ref{thm:main}}
\label{sec:out1}
 Let $K$ be a number field with an elliptic curve of rank 1.  The key to the proof of Theorem \ref{thm:main}, that is the key to the construction of a big subring of $K$ where the theorem holds, is the choice of $K$-primes to invert in the ring.  In \cite{Po2}  and \cite{PS} the inverted primes were chosen so that only a specific sequence of the elliptic curve points had its coordinates in the ring. (We remind the reader that an element of our number field is in the ring if and only if all the primes occurring in the denominator of its divisor are inverted in the ring.) In our case, almost no point of the elliptic curve will have its coordinates in the ring and we will have to represent each coordinate  by a pair consisting of a ``numerator'' and the corresponding ``denominator''.  This is the reason for having a class Diophantine model at the end instead of a regular Diophantine model: every coordinate of an elliptic curve point will be represented by an equivalence class of pairs of ``numerators'' and ``denominators'', as in a standard construction of the fraction field of a ring.

 To explain the main ideas of the proof we for the moment simplify the situation assuming that $K=\Q$, there are no torsion points,  and every non-trivial multiple of the generator $P$ has a primitive divisor.   In other words we assume that for every $n >0$, there exists a prime dividing the reduced denominators of the affine coordinates of $[n]P$ such that this prime does not divide the reduced denominators of the coordinates of any $[m]P$ with $0 <m<n$.   (In general this will be true for sufficiently large $n$ only.  See Proposition \ref{prop:biggerS}.)  We will also assume that the coordinates of $P$ itself are non-zero integers.   (In ``real life'' we will invert the primes which appear in the denominator of the coordinates of $P$. Also the primitive divisor requirement and the chosen form of the Weierstrass equation will force all the non-trivial multiples of $P$ to have non-zero coordinates.)  Under our assumptions we can represent $[n]P$ for a non-zero integer $n$, as a pair $\displaystyle \left (\frac{U_n}{V_n},\frac{X_n}{Y_n}\right )$, where $U_n \not = 0, V_n>0, X_n \not=0, Y_n >0$ are integers and $(U_n,V_n)=1, (X_n,Y_n)=1$. Later we will not be able to assume that $U_n, V_n, X_n, Y_n$ are integers but only that these are elements of our big ring.  However, we will able to treat the variables ranging over the big rings almost in the same way as if they were integers.

 From \cite{Po2}  (see Proposition  \ref{le:Po3.1}, Lemma \ref{le:orderchange}, and Lemma \ref{cor:intersec} in this paper) we know that
 \begin{equation}
 \label{divcond}
 \mbox{if $m, n \in \Z_{\not=0}$ with $m |n$, then $V_m | V_n$ in $\Z$, and conversely if $V_m | V_n$ in $\Z$ then $m | n$},
 \end{equation}
 and
  \begin{equation}
  \label{relprime}
  \mbox{if $k,m \in \Z_{>0}$, then all the primes occurring in $(V_k,V_m)$ occur in $V_{(k,m)}$},
  \end{equation}
where $(V_k,V_m)= \mbox{GCD}(V_k,V_m)$ and $(k,m)= \mbox{GCD}(k,m)$ in $\Z$. (Since we assumed that every non-trivial multiple of $P$ has a primitive divisor, we do not have to worry about $k$ and $m$ being large enough.)   Given our assumptions on the coordinates of $P$, we have that $V_1=1$, and if $(k,m)=1$, then $(V_k,V_m) =1$. Thus, if $k$ and $m$ are non-zero relatively prime integers, then $V_kV_m | V_{km}$.  Unfortunately, in general $V_{km}$ does not divide $V_kV_m$.  In particular, $V_{km}$ is divisible by some prime powers which do not occur in $V_k$ and $V_m$. {\it So the main idea behind the proof is to invert these extra primes to force $V_{km}$ to divide $V_kV_m$ in the resulting ring.}  Of course we have to leave enough primes uninverted so that \eqref{divcond} still holds in the ring.

 We now describe the primes we do not invert.  For each rational prime $p$ and any positive integer $\ell$ we keep {\it uninverted} the largest primitive divisor of $[p^{\ell}]P$. We call these primes {\it indicator primes}.  (The idea that the indicator primes are enough to identify uniquely positive multiples of a generator was first investigated in \cite{CS}.)   We invert all the other primes and denote by $R$ the resulting subring of $\Q$.  Observe that for $m=\prod p_i^{\ell_i}$, we have that $V_m$ is divisible by the indicator prime of each $[p_i^{\ell}]P$ for all $i$ and all $\ell=1,\ldots, \ell_i$, and, because of \eqref{relprime}, by no other  indicator primes.  Indeed, first suppose $q$ is an indicator prime $V_{p^{r}}$, where $p \not = p_i$ for any $i$. In this case by \eqref{relprime}, $q$ divides $V_{(p^r,m)}=V_1=1$ and we have a contradiction.  Next assume that $q$ is an indicator prime for some $V_{p_i^r}$, where $r >\ell_i$.  By definition of an indicator prime we have that $q | V_{p_i^r}$ but $q$ does not divide $V_{p_i^j}$ for any $j \in \{1, \ldots, r-1\}$.  Applying \eqref{relprime} again we obtain $q | V_{(p^r,m)}=V_{p^{\ell_i}}$ contradicting our assumptions in this case also.

 So now we are in a situation where for $k, m \in \Z_{\not =0}$ and relatively prime, $V_kV_m$ and $V_{km}$ are divisible by the same uninverted primes.  Unfortunately, there is one more point to take care of.  The indicator primes do not necessarily appear to the same power in $V_k$, $V_m$ and $V_{km}$ (see Lemma \ref{le:orderchangeold} in this paper).
Here we need another technical result from \cite{Po2}.   Let $q\not =p$ be rational primes.
  \begin{equation}
  \label{power}
  \mbox{If $\ord_{q}V_m >0$, then $\ord_qV_{pm}= \ord_qV_m$, and $\ord_qV_{qm}>\ord_qV_m $}.
  \end{equation}
   Therefore we will need another condition on $k$ and $m$ besides being relatively prime: $(k,V_m)_R=1$ and $(m,V_k)_R=1$.   With these additional assumptions we conclude that that $V_{km} |_R V_kV_m$ in our ring.  (Here for $A, B \in R$ we write ``$(A,B)_R=1$'' to  mean that the reduced numerators of $A$ and $B$ are not simultaneously divisible by any non-inverted prime, and we write ``$A |_R B$'' to indicate the divisibility in the ring, i.e. the fact that $\displaystyle \frac{B}{A} \in R$.)  To summarize the discussion above we can now say $\forall n, k, m \in \Z_{\not=0}:$
 \[
   \begin{array}{c}
   [(k,m)=1\land  (k,V_m)_R=1\land(m,V_K)_R=1] \\
 \Big  \Downarrow \\
  (V_n |_RV_mV_k \land V_mV_k|_R V_n \Leftrightarrow  |n|=|km|).
   \end{array}
\]
    (See Proposition \ref{prop:canmultiply} and Lemma \ref{le:product}.)

   If $(k,m)=1, (k,V_m)_R=1$, and $(m,V_k)_R=1$, we say that the indices $k$ and $m$ can be ``multiplied directly''.   Before we explain how to ``multiply'' arbitrary indices, note that for {\it any} triple of non-zero indices $k, m$ and $n$  we have that
 \begin{equation}
V_n |_R V_mV_k \mbox{ and } V_mV_k|_R V_n \mbox{ implies } |n|=|km|.
\end{equation}

(As above, the divisibility bar with a subscript $R$  here refers to the divisibility in our ring.)   Note also that, as a general matter, for any ring of characteristic not equal to 2, to define multiplication, it is enough to define squaring: $xy = \frac{1}{2}((x+y)^2-x^2-y^2)$.  

To take care of the indices that we cannot multiply directly we show that for every even index $k$ there exists an odd integer $w$ such that pairs $k$ and $w$ and $k$ and $k+w$ can be multiplied directly. (See Proposition \ref{prop:directly} and Remark \ref{rem:directly}.) In other words we are able to say, given an index $k \in 2\Z_{\not =0}$, that there exists a $w \in \Z_{\not =0}$, such that $GCD(k,w)=1$ and for some  $s, t \in \Z_{\not = 0}$ we have that
\begin{equation}
V_{k+w}V_w |_R V_s \land V_s |_R V_{k+w}V_w
\end{equation}
and
\begin{equation}
V_kV_w |_R V_t \land V_t |_R V_kV_w
\end{equation}
or, in other words,
\begin{equation}
\label{eq:abs}
 |(k+w)k|=|s| \land |kw|=|t|
\end{equation}
If not for absolute values in \eqref{eq:abs}, we would be done, since we would be able to define a square of $k$.    We deal with absolute values via considering all possible cases and using \eqref{divcond} in Lemma \ref{le:define}.

Over $\Z$, given integers  $U>0, V>0, X\not =0, Y\not = 0$  such  that $\displaystyle \frac{U}{V}$ and $\displaystyle \frac{X}{Y}$ satisfy the chosen Weierstrass equation and such that $(U,V)=1, (X,Y)=1$, we can conclude that  $(U,V,X,Y)=(U_n, V_n, X_n, Y_n)$ for some unique $n \in \Z_{\not=0}$.  Unfortunately, if we now assume that $(U>0, V>0, X\not =0, Y\not = 0) \in R^4$, where $R$ is, as above, our ring with infinitely many primes inverted, and $\displaystyle \frac{U}{V}$ and $\displaystyle \frac{X}{Y}$ satisfy the chosen Weierstrass equation with $$(U,V)_R=1, (X,Y)_R=1,$$ then we will be able to conclude only that $U=\tilde U_n=U_n\bar U_n, V=\tilde V_n=V_n\bar V_n$, where $\bar U_n, \bar V_n$ are rational numbers whose reduced numerators and  denominators are divisible by the inverted primes only.  (A similar conclusion will apply to $(X,Y)$.)  However, since we are only interested in the divisibility by the non-inverted indicator primes, the ``bar'' parts do not matter  or in other words, for any $k, m, n \in \Z_{\not = 0}$ we still have that
 \[
   (k,m)=1, (k,V_m)=1, (m,V_k)=1 \Longrightarrow   (\tilde V_n |_R\tilde V_m\tilde V_k \land \tilde V_m\tilde V_k|_R \tilde V_n \Leftrightarrow  |n|=|km|).
   \]
 This is so, because $(k,m)=1, (k,V_m)=1, (m,V_k)=1 \Longleftrightarrow (\tilde V_k, \tilde V_m)_R=1$ and $\tilde V_n |_R\tilde V_m\tilde V_k \Longleftrightarrow  V_n |_R V_mV_k$, etc.   The fact that we can express the condition of being relatively prime in our ring in polynomial terms is demonstrated in Lemma \ref{le:reprime}.  Unfortunately, when the underlying field has a class number greater than one, there are other technical complications requiring raising variables to the power divisible by a class number to obtain relatively prime numerators and denominators.  (See Notation \ref{not:2}, Item \ref{it:d} and Remark \ref{rem:d}.)

    The last point that needs to be explained is the density of the inverted and the non-inverted prime sets.  In \cite{Po2} and \cite{PS}, it was shown that the natural density of the indicator primes corresponding to the prime multiples of any infinite order point is 0.  So the only remaining question is the density of the indicator primes corresponding to prime power multiples of such a point, when the power is at least 2.  This density is also 0 and the corresponding calculation is  much easier.  It was first carried out in \cite{CS} and is reproduced in the appendix of this paper for the convenience of the reader.

\section{An Outline of the Proof of Theorem \ref{thm:first-order}.}
In this section we keep for the moment the simplifying assumptions and notation of the preceding section, i.e. we assume that we are dealing with a rank one elliptic curve over $\Q$ with a trivial torsion group, and a Weierstrass equation as above, and that every non-trivial multiple of a generator has a primitive divisor.  We also assume that Theorem \ref{thm:main} holds or in other words in a big subring $R$ of $\Q$ described above we have defined existentially multiplication of indices.

If $x \in \Q$ and $\displaystyle x = \frac{A}{B}$, where $A, B \in R$ with $AB \not = 0$, then we say that $A$ and $B$ are a reduced numerator and a denominator respectively, if $(A,B)_R= 1$.  In other words, neither $A$, nor $B$ are divisible by ``extra'' non-inverted primes.  If $R=\Z$, this definition is the same as the usual one.  We now need the following results from \cite{Po} (Lemmas \ref{le:anydivisor} and \ref{le:equiv} of this paper):
\begin{equation}
\label{eq:divide3}
\begin{array}{c}
\mbox{ ``For any sufficiently large }l \in \Z_{>0}, \mbox{ for any }k \in \Z_{>0} \mbox{ we have that }\\
\mbox{the reduced denominator of } x_l \mbox{ divides the reduced numerator of }\left (\displaystyle \frac{x_{l}}{x_{kl}}-k^2\right)^2 \mbox{ in } \Z\mbox{''},
\end{array}
\end{equation}and
\begin{equation}
\label{eq:divide4}
\begin{array}{c}
\mbox{``For any } n  \in \Z_{>0} \mbox{ there exists } l \in \Z_{>0} \mbox{ such that } \\
 n \mbox{ divides the reduced denominator of } x_l \mbox{ in } \Z\mbox{''}.
\end{array}
\end{equation}
 Now let $z$ be an arbitrary element of our big ring with the following property:   there exists   a non-zero integer $k$,  such that  for all rational numbers $b$ in our ring,  there exist non-zero integers $i$ and $j$ satisfying the equations \eqref{eq:divide2}--\eqref{eq:congruence1} below.
\begin{equation}
\label{eq:divide2}
b^2 \mbox{ divides the reduced denominator of } x_i \mbox{ in our ring.}
\end{equation}%
\begin{equation}
\label{eq:multiply1}
j=ik
\end{equation}%
\begin{equation}%
\label{eq:congruence1}
\begin{array}{c}
 \mbox{ The reduced denominator of } x_i\\
  \mbox{ divides the reduced numerator of } \displaystyle (z - \frac{x_i}{x_j})^2 \mbox{ in our ring.}
  \end{array}
\end{equation}
(Here, as above, $x_k, x_i, x_j$ are the $x$-coordinates of $[k]P, [i]P$ and $[j]P$ respectively.)  Then $z \in \Z$.\\

Conversely, if  $z$ above is a square of a non-zero integer, then we can find a $k \in \Z_{\not = 0}$ such that for every $b$ in our big ring there exist $i$ and $j$ so that \eqref{eq:divide2} -- \eqref{eq:congruence1} are satisfied.\\

First assume that $z$, a rational  number in our ring, is fixed.  Let $k$ be the corresponding non-zero integer, $b$ an arbitrary element of the ring and assume that  $i,j \in \Z_{\not = 0}$ are such that  the equations above are satisfied.  From \eqref{eq:divide3} and \eqref{eq:divide2} we conclude that $b$ divides the reduced numerator of  $\displaystyle \left (\frac{x_{i}}{x_{j}}-k^2\right) = \left (\frac{x_{i}}{x_{ik}}-k^2\right)$ as well as the reduced numerator of  $\displaystyle (z - \frac{x_i}{x_j})=(z - \frac{x_i}{x_{ik}})$ in our ring.  Thus, $b$ divides the reduced numerator of $z-k^2$ in our ring.  If $\displaystyle z=\frac{z_1}{z_2}$, where $z_1, z_2 \in \Z_{\not = 0}$, then $b$ divides $z_1-z_2k^2$ in our ring. If we  pick $b$ to be divisible by $q^m$, where $q$ is a prime which is not inverted in our ring and $m$ is a positive integer large enough so that $q^m >|z_1-z_2k^2|$, then $q^m$ divides $z_1-z_2k^2$ in $\Z$ and the only way the divisibility condition can hold is for $z_1=z_2k$.  Without loss of generality we can assume that $z_1$ and $z_2$ were picked to be relatively prime in $\Z$, and since $k$ is a non-zero integer, we must conclude that $z_2=1$, and $z=z_1=k^2$.

Assume now that $z=k^2$ where $k \in \Z_{\not =0}$.  Let $b$ be any rational number in our ring.  Let $i>0$ be such that $b^2$ divides the reduced denominator of $x_i$ and $i$ is sufficiently large so that \eqref{eq:divide3} holds for $l=i$.  Such an $i$ exists by \eqref{eq:divide4}.  Finally let $j=ik$ and observe that \eqref{eq:congruence1} now holds by \eqref{eq:divide3}.

\section{Elliptic curves}
\setcounter{equation}{0}
We now proceed with the detailed description of the proof.  In this section we lay down the elliptic curve foundations of our results.  Many  of the technical
details in this section  are taken from \cite{CS}, \cite{Po}, \cite{Po2} and \cite{PS}.  Below we indicate which technical results have been taken from other papers.
\begin{notationassumptions}%
\label{S:notation section}%
We add the following notation and assumptions to the list above.
\begin{itemize}%
\item Let $E$ be an elliptic curve of rank~1 defined over $K$ (in particular, we assume such an $E$ exists).%
\item We fix a Weierstrass equation $W: y^2=x^3+ax+b$ for $E$ with all the coefficients in the ring of integers of $K$. %
\item Let $E(K)_\tors$ be the torsion subgroup of $E(K)$.%
\item Let $t$ be a multiple of $\#E(K)_\tors$.%
\item Let $Q \in E(K)$ be such that $Q$ generates $E(K)/E(K)_\tors$.%
\item Let $P:=[t]Q$.%

\item Let ${\mathcal S}_\bad = \calS_{\bad}(W,P,K) \subseteq \calP_K$ consist of the primes that ramify in $K/\Q$,
the primes for which the reduction of the chosen Weierstrass model is singular (this includes all primes above
$2$), and the primes at which the coordinates of $P$ are not integral.%
\item For $n \in \Z_{\not =0}$ write $[n]P=(x_n,y_n)=(x_n(P), y_n(P))$ where $x_n,y_n \in K$.%
\item For $n \in \Z_{\not = 0}$, let the divisor of $x_n(P)$ be of the form%
\[%
\frac{{\mathfrak a}_n}{{\mathfrak d}_n}{\mathfrak b}_n=\frac{{\mathfrak a}_n(P)}{{\mathfrak
d}_n(P)}{\mathfrak b}_n(P)
\]%
where
\begin{itemize}%
\item ${\mathfrak d}_n=\prod_{{\mathfrak q}}{\mathfrak q}^{-a_{{\mathfrak q}}}$, where the product
is taken over all primes ${\mathfrak q}$ of $K$ not in $\calS_{\bad}$ such that $a_{{\mathfrak
q}}=\ord_{{\mathfrak q}}x_n <0$.%
\item ${\mathfrak a}_n=\prod_{{\mathfrak q}}{\mathfrak q}^{a_{{\mathfrak q}}}$, where the product is
taken over all primes ${\mathfrak q}$ of $K$ not in $\calS_{\bad}$ such that $a_{{\mathfrak
q}}=\ord_{{\mathfrak q}}x_n >0$.%
\item ${\mathfrak b}_n = \prod_{{\mathfrak q}}{\mathfrak q}^{a_{{\mathfrak q}}}$, where the product
is taken over all primes ${\mathfrak q}\in \calS_{\bad}$ and $a_{{\mathfrak q}}=\ord_{{\mathfrak
q}}x_n$.%
\end{itemize}%
\item For $n$ as above, let ${\calS}_n = {\calS}_n(P)=\{\pp \in {\mathcal P}_K : \pp | {\mathfrak
d}_n\}$. By
definition of $\calS_\bad$ and ${\mathfrak d}_n$, we have $\calS_1=\emptyset$.%
\item For $\ell \in \calP_{\Q}$, define $a_\ell$ to be the smallest positive integer such that for any $j \geq a_{\ell}$ we have that
$\calS_{\ell^{j}} \setminus \calS_{\ell^{j-1}}\not = \emptyset$.  By  Proposition \ref{prop:biggerS} below, for all
but finitely many primes $\ell$ we have that $a_\ell =1$ and for all $\ell$ we have that $a_{\ell}$ is well defined.
\item For $j \in \Z_{\geq 1}$, let $\pp_{\ell^{j}}(P)=\pp_{\ell^j}$ be a prime of the largest norm
in $\calS_{\ell^{j}} \setminus \calS_{\ell^{j-1}}$, if such a prime exists. (This prime will be called the indicator prime for $[\pp_{\ell^{j}}]P$.)
\item Let $m_0 =\prod_{a_{\ell} >1}\ell^{a_{\ell}-1}$. (Note that $m_0$ is well defined since, as we have observed above, for all
but finitely many primes $\ell$ we have that $a_\ell =1$.)
\item For all $j \in \Z_{\geq 1}$ let $\qq_{\ell^j}=\pp_{\ell^{j + \ord_{\ell}m_0}}$.
\item Let $T=[m_0]P$.
\item Let $\calV_K  =\calV_K(P) = \{\pp_{\ell^j}: \ell \in \calP_{\Q}, j \in \Z_{>0}\}$.%
\item Let $\calW_K =(\calP_K \setminus \calV_K) \cup \calS_{m_0}$.  ($\calW_K$ will be the set of the inverted primes.)
\item Let $\calC_n = (\calS_n \cap \calV_K) \setminus \calS_{m_0}$.   Note that $\calC_{m_0} = \emptyset$.  ($\calC_n$ will be  the collection of the prime factors of the ``$|n|$''-th denominator which are not inverted.)
\item Let $\calX_n=\calS_{m_0n}$. ($\calX_n$  will be the set of the ``not-bad denominator primes'' for $[n]T$.)
\item Let $\calY_n = \calC_{m_0n}$ and observe that $\calY_1$ is empty.   ($\calY_n$ will be the set of the non-inverted ``denominator'' primes for $[n]T$.)
\item Let ${\mathfrak c}_n=\prod_{{\mathfrak q}}{\mathfrak q}^{-a_{{\mathfrak q}}}$, where the product is taken over all primes ${\mathfrak q}$ of $K$ not in $\calW_K$ such that $a_{{\mathfrak q}}=\ord_{{\mathfrak q}}x_n <0$. (The divisor $\cc_n$ will be the non-inverted part of the  ``$|n|$''-th denominator.)
 \item Let ${\mathfrak f}_n={\mathfrak c}_{m_0n}$.
 \item For $x \in K$, let $\dd(x)= \prod_{{\mathfrak q}}{\mathfrak q}^{-a_{{\mathfrak q}}}$, where the product is taken over all primes ${\mathfrak q}$ of $K$ such that $a_{{\mathfrak q}}=\ord_{{\mathfrak q}}x <0$.  Let $\nn(x) = \dd(x^{-1})$.
\item For $x \in K$, let $\dd_{\calW_K}(x)= \prod_{{\mathfrak q}}{\mathfrak q}^{-a_{{\mathfrak q}}}$, where the product is taken over all primes ${\mathfrak q}$ of $K$ not in $\calW_K$ such that $a_{{\mathfrak q}}=\ord_{{\mathfrak q}}x <0$.  Let $\nn_{\calW_K}(x) = \dd_{\calW_K}(x^{-1})$.
\end{itemize}%
\end{notationassumptions}%
Below we combine ideas from \cite{Po}, \cite{PS} and \cite{CS} to show that it is enough to have one non-inverted indicator prime for every prime power of the index to identify the index of a point uniquely (up to a sign).  At the same time, if we don't invert only the indicator primes of the index prime powers, we will have ``almost'' arranged for the multiplication of indices.

As pointed out above, denominator prime sets are not enough to establish a sign of an index.  This is demonstrated by the lemma below.
\begin{lemma}%
\label{le:minus}%
For any $n \in \Z_{\not=0}$ we have that $\calS_n=\calS_{-n}$,$\calC_n = \calC_{-n}, \calX_{-n}=\calX_{n}$,  $\calY_{-n}=\calY_{n}$, and $\ff_n =\ff_{-n}$.
\end{lemma}%
\begin{proof}%
Given the choice of our Weierstrass equation, we have that $x_{-n}= x_n$.
\end{proof}%
Our next step is to establish several important properties of the primes which appear in the denominators in Propositions \ref{le:Po3.1}--\ref{cor:equality2}.  Fortunately for us, most of the technical work has already been done elsewhere.
\begin{proposition}[Lemma 3.1 of \cite{PS}]%
\label{le:Po3.1}%
Let ${\mathfrak R}$ be an integral divisor of $K$. Then $$\{n \in \Z \setminus \{0\}: {\mathfrak R}
\mid {\mathfrak d}_n(P)\} \cup \{0\}$$ is a subgroup of $\Z$.
\end{proposition}%

\begin{proposition}[Proposition 3.5 of \cite{PS}]
\label{prop:biggerS}
There exists $C >0$ such that for all $\ell, m \in \calP_{\Q}$ with $\max(\ell,m) >C$ we have that
$\calS_{\ell m} \setminus (\calS_\ell \cup \calS_m) \ne \emptyset$.
\end{proposition}

\begin{lemma}
\label{le:orderchangeold}
Let $n \in \Z_{\ge 1}$.
Suppose that ${\mathfrak t} \in \calP_K$ divides ${\mathfrak d}_n$,
and $p \ge 2$ is a rational prime.
\be
 \item If ${\mathfrak t} \mid p$, then $\ord_{{\mathfrak t}}{\mathfrak d}_{pn} \geq 2+ \ord_{{\mathfrak t}}{\mathfrak d}_n$.
 \item If ${\mathfrak t} \nmid p$, then $\ord_{{\mathfrak t}}{\mathfrak d}_{pn}= \ord_{{\mathfrak t}}{\mathfrak d}_n$.
\ee
\end{lemma}
\begin{proof}%
The proof of the lemma is almost identical to the proof of Lemma 3.3 of \cite{PS} except for the fact that we allow $p=2$.  We also remind the reader that any $\ttt$ dividing $\dd_n$ is automatically not in $\calS_{\bad}$ and therefore is not dyadic, ramified over $\Q$ or is among primes at which our Weierstrass model has a bad reduction.
\end{proof}%
\begin{corollary}
  \label{le:orderchange}
Let $n \in \Z_{\ge 1}$.
Suppose that ${\mathfrak t} \in \calP_K$ divides ${\mathfrak c}_n$  (or $\ff_n$), and $p \ge 2$ is a rational prime.
\be
 \item If ${\mathfrak t} \mid p$, then $\ord_{{\mathfrak t}}{\mathfrak c}_{pn}\geq  2+ \ord_{{\mathfrak t}}{\mathfrak c}_n$ (or $\ord_{{\mathfrak t}}{\mathfrak f}_{pn} \geq2+\ord_{{\mathfrak t}}{\mathfrak f}_n$).
 \item If ${\mathfrak t} \nmid p$, then $\ord_{{\mathfrak t}}{\mathfrak c}_{pn}= \ord_{{\mathfrak t}}{\mathfrak c}_n$ (or $\ord_{{\mathfrak t}}{\mathfrak f}_{pn}= \ord_{{\mathfrak t}}{\mathfrak f}_n$).
\ee
\end{corollary}
\begin{proof}
The corollary follows immediately from the lemma above if  we note that we obtain $\cc_n$ from $\dd_n$ by removing factors of $\dd_n$ which are in $\calW_K$, and $\ff_{n}=\cc_{m_0n}$.
\end{proof}%

\begin{lemma}[Lemma 10 of \cite{Po}]%
\label{le:anydivisor}%
Let $\mathfrak A$ be any integral divisor of $K$.  Then there exists  $k \in \Z>0$ such that $\mathfrak A \Big{|} \dd(x_k)$.
\end{lemma}%

\begin{lemma}%
\label{cor:intersec}
Let $m,n \in \Z \setminus \{0\}$, and let $(m,n)$ be their GCD. Then
\[
{\mathcal S}_m \cap {\mathcal S}_n ={\mathcal S}_{(m,n)},
\]
\[
\calX_m \cap \calX_n = \calX_{(m,n)},
\]
\[
\calC_m \cap \calC_n = \calC_{(m,n)},
\]
and
\[
\calY_m \cap \calY_n = \calY_{(m,n)}.
\]
In particular, if $(m,n)=1$, then
\[
{\mathcal S}_m \cap {\mathcal S}_n=\emptyset,
\]
\[
{\mathcal X}_m \cap {\mathcal X}_n=\calX_1=\calS_{m_0},
\]
\[
\calC_m \cap \calC_n = \calC_{1}=\emptyset,
\]
and
\[
{\mathcal Y}_m \cap {\mathcal Y}_n=\calC_{m_0}=\emptyset.
\]
\end{lemma}%
\begin{proof}
The assertion ${\mathcal S}_m \cap {\mathcal S}_n = {\mathcal S}_{(m,n)}$ is exactly Lemma 3.2 of \cite{PS}.  Therefore if $(m,n)=1$ we have that $\calS_{(m,n)}=\calS_1=\emptyset$ by definition of $\calS_n$.  Further, by definition,
\[
\calX_n=\calS_{m_0n}, \calX_m=\calS_{m_0n}
\]
 and therefore,
\[
\calX_m \cap \calX_n=\calS_{m_0n} \cap \calS_{m_0m}=\calS_{m_0(m,n)}=\calX_{(m,n)}.
\]
Thus, if $(m,n)=1$ we have
\[
\calX_m \cap \calX_n=\calX_1=\calS_{m_0}.
\]
Also by definition,
\[
\calC_n=(\calS_n \cap \calV_K) \setminus \calS_{m_0}, \calC_m=(\calS_m \cap \calV_K)\setminus \calS_{m_0}
\]
and therefore,
\[
\calC_n \cap \calC_m = (\calS_m \cap \calS_n \cap \calV_K)\setminus \calS_{m_0}= (\calS_{(m,n)} \cap \calV_K)\setminus \calS_{m_0}=\calC_{(m,n)}.
\]
Consequently, if $(m,n)=1$ we have that
\[
\calC_m \cap \calC_n=\calC_1 =(\calS_1 \cap \calV_K)\setminus \calS_{m_0}= \emptyset.
\]
Finally, again by definition,
\[
\calY_n= \calC_{m_0n}, \calY_m= \calC_{m_0m}
\]
and therefore,
\[
\calY_n \cap \calY_m = \calC_{m_0n}\cap \calC_{m_0m}=\calC_{m_0(m,n)}=\calY_{(m,n)}.
\]
Consequently, if $(m,n)=1$ we have that
\[
\calY_m \cap \calY_n=\calY_1=\calC_{m_0}=(\calS_{m_0} \cap \calV_K)\setminus \calS_{m_0}=\emptyset.
\]

\end{proof}

\begin{corollary}%
\label{cor:div}%
For any $\ell \in \calP(\Q)$ and any $j \in \Z_{>0}$ we have that $\qq_{\ell^j}$ exists, and
$\displaystyle \qq_{\ell^j} \in \calX_k=\calS_{km_0}$ if and only if $\ell^j$ divides $k$.  (We remind the reader that by definition, $\qq_{\ell^j}=\pp_{\ell^{j +\ord_{\ell}(m_0)}}$ is the indicator prime of $[\ell^{j + \ord_{\ell}m_0}]P$.)
\end{corollary}%
\begin{proof}%
By definition of $\qq_{\ell^j}$, to establish its existence it is enough to show that $$\calS_{\ell^{\ord_{\ell}m_0+j}}\setminus \calS_{\ell^{j-1+\ord_{\ell}m_0}} \not = \emptyset.$$ At the same time, from the definitions of $m_0$ and $a_{\ell}$ we have that
\[%
 \calS_{\ell^{\ord_{\ell}m_0+j}}\setminus \calS_{\ell^{j-1+\ord_{\ell}m_0}}= \calS_{\ell^{a_{\ell}-1+j}}\setminus \calS_{\ell^{j-1+a_{\ell}-1}}\not =\emptyset,
\]%
 and therefore $\qq_{\ell^j}$ exists.

 Now suppose $j >0$ and  $\qq_{\ell^j} \in \calX_k=\calS_{km_0}$. Then by definition of $\qq_{\ell^j}$, we have that
\[%
\pp_{\ell^{j+\ord_{\ell}m_0}}\in \calS_{km_0} \cap \calS_{\ell^{j+\ord_{\ell}m_0}}=\calS_{GCD(km_0,
\ell^{j+\ord_{\ell}m_0})} \subseteq \calS_{\ell^{\ord_{\ell}(km_0)}}
\]%
by Lemma \ref{cor:intersec}. But by the same lemma and the definition of indicator primes, $$ \pp_{\ell^{j+\ord_{\ell}m_0}} \in \calS_{\ell^{\ord_{\ell}(km_0)}}\Leftrightarrow j \leq \ord_{\ell}k.$$

Conversely,  suppose $j >0$ and $j \leq \ord_{\ell}k$.  Then $\pp_{\ell^{j+\ord_{\ell}m_0}}\in \calS_{\ell^{j+\ord_{\ell}m_0}}
\subset \calS_{km_0}$ by Lemma \ref{cor:intersec} once again and $\qq_{\ell^j} \in \calX_k$.
\end{proof}%

 \begin{corollary}%
 \label{cor:div2}  %
 \be%
\item For any $k \in \Z_{>1}$ we have that
 \[%
 \calY_k = \{ \qq_{\ell^j}: \ell \in \calP_{\Q}, 0< j \leq \ord_{\ell}k \}.
 \]%
\item For $k,n \in \Z_{>1}$ we have that $\calY_k \subseteq \calY_n$ if and only if $k \Big{|} n$.
\item  For $k, n \in \Z_{>1}$ we have that $\ff_k \Big{|} \ff_n$ if and only if $k   \Big{|} n$.%
\item For $k, n \in \Z_{>1}$ we have that $(k,n)=1$ if and only if $(\ff_k,\ff_n)=(1)$, where $(1)$ is a trivial divisor.
\ee%
 \end{corollary}
 \begin{proof}%
 \be%
\item First we observe that by definition of  $\calY_k=\calC_{m_0k} = \calS_{m_0k} \setminus \calW_K=\calX_k \setminus \calW_K$,  these prime sets contain only the primes of
the form $\pp_{\ell^j}$ for some $\ell \in \calP_{\Q}$ and some $j \in \Z_{>0}$.  Secondly, by Corollary \ref{cor:div}, we also have that $\qq_{\ell^j}  \in \calX_k$ if and only if  $0< j \leq \ord_{\ell}k$.

\item If we assume that $k \Big{|} n$, then $ \calX_k \subseteq \calX_n$ by Lemma \ref{cor:intersec} and consequently, $\calY_k \subseteq \calY_n$.
Conversely, if we suppose that $\calY_k \subset \calY_n$, then for every rational prime $\ell$ we have that $\qq_{\ell^{\ord_{\ell}(k) }} \in \calY_{n}$ by  Part 1 of this corollary. Thus, by  Part 1 again,  for
every rational prime $\ell$ we have that $ \ell^{\ord_{\ell}(k) } \Big{|} n$. Consequently $k$ divides $n$.%
 \item If we first assume that   $\ff_{k} \Big{|} \ff_{n}$,  then  $\calY_{k} \subseteq \calY_{n}$ and $k \Big{|} n$ by Part 2 of this corollary.  Next if we suppose $k \Big{|} n$,  then $\calY_{k} \subseteq \calY_{n}$  by Part 2 of this corollary again,  and consequently  $\ff_{k} \Big{|} \ff_{n}$   by Corollary \ref{le:orderchange}.
 \item Suppose $(k,n) = 1$, then $\calY_k \cap \calY_n = \emptyset$ by Corollary \ref{cor:intersec}. Since all the prime divisors of $\ff_k$ are in $\calY_k$, and all the prime divisors of $\ff_n$ are in $\calY_n$, we must conclude that $(\ff_k,\ff_n)=(1)$.  Conversely, if $(\ff_k,\ff_n)=(1)$, then $\calY_k \cap \calY_n = \emptyset=\calY_{(k,n)}$, where the last equality holds by Corollary \ref{cor:intersec}.  But from Part 1, we conclude that $(k,n)=1$ since $\calY_1$ is the only $\calY_m$ with $m >0$ which is an empty set.
 \ee%
 \end{proof}%

The next corollary is the first step towards the existential definition of multiplication of indices.
\begin{corollary}%
\label{prop:factors}
Let $m, k \in \Z_{\not =0}$ with $(m,k)=1$.  Then $\calY_{mk} =\calY_m \cup \calY_k$
\end{corollary}%
\begin{proof}%
Since $(m,k)=1$ the assertion  follows from the Part 1 of Corollary \ref{cor:div2}. Indeed, for any $j \in \Z_{>0}$ and $\ell \in \calP_{\Q}$ we have that $0 < j \leq \ord_{\ell}mk$ if and only if either  $0 < j \leq \ord_{\ell}m$ or   $0 < j \leq \ord_{\ell}k$.
\end{proof}%
While we established already that the denominator prime sets cannot distinguish between positive and negative indices, the result below tells us that the indicator primes identify the {\it absolute value} of the index for a multiple of $T$ uniquely.
\begin{corollary}%
\label{cor:equality1}
Let $n_1, n_2 \in \Z_{>0}$ be such that    $\calY_{n_1}= \calY_{n_2}$.  Then $n_1=n_2$,
\end{corollary}%
\begin{proof}%
By Corollary \ref{cor:div2} we have that $n_1$ divides $n_2$ and $n_2$ divides $n_1$. Thus, $n_1=n_2$.
\end{proof}%
From Corollary \ref{cor:equality1} we immediately obtain the proposition below.
\begin{corollary}%
\label{cor:equality2}
Let $n_1, n_2 \in \Z_{>0}$ be such that    $\ff_{n_1}= \ff_{n_2}$.  Then $n_1=n_2$,
\end{corollary}%
\begin{proof}%
The equality $\ff_{n_1}= \ff_{n_2}$ implies   $\calY_{n_1}= \calY_{n_2}$ and we are done by Corollary \ref{cor:equality1}.
\end{proof}%
We are now ready to conclude that under our definitions and under certain relative primality assumptions, the denominator of the product is equal to the product of the denominators.
\begin{proposition}%
\label{prop:canmultiply}
 Let $m, k \in \Z_{>0}$  with $(m,k)=1$ and $((m), \ff_{k}) = 1, ((k),\ff_{m})=1$.  Then $\ff_{mk} = \ff_{k}\ff_{m}$. (Here we consider $(k), (m)$ as divisors in $K$.)
\end{proposition}   %
\begin{proof}%
By Corollary \ref{le:orderchange} and Corollary \ref{prop:factors} we have that $\ff_{k}\ff_{m}$ divides $\ff_{mk}$.
Thus, it is enough to show that $\ff_{mk}$ divides $\ff_{k}\ff_{m}$. So let $\pp \in \calP_K$ be such that
$\ord_{\pp}\ff_{mk}=a>0$. Then by Corollary \ref{prop:factors} either $\ord_{\pp}\ff_{m} >0$ or $\ord_{\pp}\ff_{k} >0$
but both inequalities cannot hold at the same time since $(k,m)=1$. (See Lemma \ref{cor:intersec}.) Without loss of generality, assume the first alternative
holds. By assumption $\ord_{\pp}k = 0$ and therefore by Corollary \ref{le:orderchange} we  have that $\ord_{\pp}\ff_{mk}
=\ord_{\pp}\ff_{m}$.
 \end{proof}%
 \begin{definition}%
 Let $m, k \in \Z_{>0}$ be such that  $(m,k)=1$ and $((m), \ff_k) = 1, ((k),\ff_m)=1$.
Then we will say that $m$ and $k$ {\it can be multiplied directly}.
 \end{definition}%
 The next lemma is a converse of sorts to the Proposition \ref{prop:canmultiply}.
 \begin{lemma}%
 \label{le:product}
 Let $m,k,n \in \Z_{>0}$, $(\ff_k,\ff_m) =1$  and $\ff_m \ff_k = \ff_n$. Then $(k,m)=1, n = mk$, and $(\ff_k,(m))= (\ff_m,(k))=1$.
 \end{lemma}%
 \begin{proof}%
 First  we show that     $(k,m)=1$.   Suppose not.  Let $\ell$ divide $(m,k)$ and conclude that $$\qq_{\ell} \in \calY_m\cap \calY_k=\calY_{(m,k)}$$ by Corollary \ref{cor:div}, and therefore $(\ff_k,\ff_m) \not =1$.  Thus $(k,m)=1$, and by assumption and Proposition \ref{prop:factors} we now have that $\calY_n=\calY_k\cup \calY_m = \calY_{mk}$.  By Corollary \ref{cor:equality1} we conclude that $n=mk$.  Suppose now without loss of generality    $(\ff_k,(m))\not = 1$.  Then for some $\pp \in \calP_K$ dividing $\ff_k$ it is the case that $\ord_{\pp}m>0$.  In this case by Corollary  \ref{le:orderchange} we have that
$\ord_{\pp}\ff_{km} >  \ord_{\pp}\ff_{k}=   \ord_{\pp}\ff_k  \ff_m$.
 \end{proof}%
We now show that it is not hard to find pairs of indices which can be multiplied directly. We start with a preliminary lemma.
\begin{lemma}
\label{le:getready}
Let $k, w$ be positive integers.  Let $\ttt_{1}, \ldots, \ttt_{m}$ be all the $K$-factors of $(k)$ not in $\calW_K$, and let $\ell_{1},\ldots, \ell_m \in \calP_{\Q}$ be such that $\ttt_{i}=\qq_{\ell_{i}^j}$ for some $j \in \Z_{>0}$.  (In other words, $\ttt_{i}$ is the indicator prime for $[\ell_{i}^j]T$.)  Assume further that
 $(w,\ell_{i})=1$ for all $i=1,\ldots,m$.  In this case $(\ff_w,(k))=(1)$.  (Here we consider $(k)$ as a divisor of $K$ as above.)
\end{lemma}
\begin{proof}
Suppose $(\ff_w,(k)) \not=1$.  In this case $\calY_w$ contains $\ttt_i=\qq_{\ell_i^j}$ for some $i=1,\ldots, m$ and some $j \in \Z_{>0}$.  However, by Corollary \ref{cor:div2} we must then conclude that $\ord_{\ell_i}w \geq j >0$ contradicting our assumptions.
\end{proof}
 \begin{proposition}%
 \label{prop:directly}
For any  $k \in 2\Z_{\not=0}$  there exists $v \in \Z_{\not=0}$ such that for  $w=kv+1$ the following conditions are satisfied:
\be%
\item \label{it:1} $(k,w) = 1$%
\item \label{it:2} $(k, k+w)=1$%
\item $((k),\ff_w)=1$%
 \item $(\ff_k, (w))=1$%
 \item $(\ff_{k+w}, (k))=1$%
\item $(\ff_k, (k+w))=1$%
 \ee %
 \end{proposition}%
 \begin{proof}%
Clearly Conditions (\ref{it:1}) and (\ref{it:2}) are satisfied by construction for any $v \in \Z_{\not=0}$.    Next let $\ttt_{1}, \ldots, \ttt_{m}$ be all the $K$-factors of $(k)$ not in $\calW_K $ let $\ell_{i} \in \calP_{\Q}$ be such that $\ttt_{i}=\qq_{\ell_{i}^j}$ for some $j \in \Z_{>0}$.  (In other words, as above $\ttt_{i}$ is the indicator prime for $[\ell_{i}^j]T$.)      Let $\ff_k =\prod \hh_u^{b_u}$ be the $K$ - prime factorization  of $\ff_k$, and for every $u$ let  $h_u$ be the rational prime below $\hh_u$.  We now rewrite the remaining conditions in terms of $v$, $\ell_{i}$, and $h_u$. It will be enough to arrange that the following conditions are satisfied for all  $ \ell_{i}, h_u$:

\be
\item $vk + 1 \not \equiv 0 \mod \ell_{i} $ (making sure that $\ff_w$ has no factors in common with $(k)$ by Lemma \ref{le:getready})
\item $vk + 1 \not \equiv 0 \mod h_u $   (making sure $w$ has no factors in common with $\ff_k$)
\item $k + vk + 1 \not \equiv 0 \mod \ell_{i}  $ (making sure $\ff_{k+w}$ has no factors in common with $(k)$ again by Lemma \ref{le:getready}. )
\item $k + vk + 1  \not \equiv 0 \mod h_u$ (making sure $k+w$ has no factors in common with $\ff_k$)
\ee

Note that for all $h_u$ and $\ell_{i}$ dividing $k$ all the conditions are automatically are satisfied. Thus, for any even
$\ell_{i}$ the conditions are satisfied. (No $h_u$ can be even by assumption on $\calW_K$.) Hence without loss of generality we
can assume that $k$ is not divisible by any $h_u$ or any $\ell_{i}$ and no $\ell_{i}$ is even. Note also that 99the
equivalences are the same across all $\ell_{i}$'s and $h_u$'s.  So repetition  of primes is not a problem. Let $g=g(h_u)$ or $g =
g(\ell_{i})$ be such that $g \not \equiv 0 $ and $g \not \equiv -k$ modulo the relevant prime. Such a $g$ exists for every
prime because all the primes are not even and so the residue fields contain at least three elements. (If $\ell_{i} =h_u=h_u'$  then the corresponding $g$'s are selected to be the same.)  Since we have assumed $k$ is not divisible by any of the primes $h_u$ or $\ell_{i}$, we can solve the congruence $vk+1\equiv g$ modulo each of the primes and use the Chinese Remainder Theorem to get a solution modulo all the primes simultaneously.
 \end{proof}%

 \begin{remark}%
 \label{rem:directly}
From Proposition  \ref{prop:directly} we conclude that for every $k \in 2\Z_{>0}$ there exists an odd $w
 \in \Z_{>0}$ such that $k$ and $w$ and $k$ and $k+w$ can be multiplied directly.
 \end{remark}   %

The remaining Propositions \ref{le:equiv} -- \ref{le:square} of this section will be necessary for defining integers using just one universal quantifier.  We start with a lemma which allows us to generate integers.

\begin{lemma}[Lemma 11 of \cite{Po}]%
\label{le:equiv}
There exists a positive integer $m_1$ such that for any positive integers $l, k $,
\begin{equation}
\label{eq:old}
\dd(x_{lm_1}) \Big{|} \nn\left (\frac{x_{lm_1}}{x_{klm_1}}-k^2\right)^2
\end{equation}
in the integral divisor semigroup of $K$.
\end{lemma}%
\begin{remark}
If we restrict our attention to the non-inverted primes only, we can rewrite \eqref{eq:old} as
\begin{equation}
\label{eq:new}
\dd_{\calW_K}(x_{lm_1}) \Big{|} \nn_{\calW_K}\left (\frac{x_{lm_1}}{x_{klm_1}}-k^2\right)^2
\end{equation}
\end{remark}
\begin{lemma}%
\label{le:relprime} With $m_1$ as in Lemma \ref{le:equiv},
$(\dd_{\calW_K}(x_{lm_1}), \nn_{\calW_K}(x_{klm_1}))=(1)$ in the integral divisor semigroup of $K$.
\end{lemma}%
\begin{proof}%
From  Lemma \ref{le:orderchangeold} and  Lemma \ref{cor:intersec} it follows that $\dd_{\calW_K}(x_{lm_1})$ divides $\dd_{\calW_K}(x_{klm_1})$ and by definition $(\dd_{\calW_K}(x_{klm_1}), \nn_{\calW_K}(x_{klm_1}))=(1)$.
\end{proof}%

From Lemma \ref{le:equiv} and Lemma \ref{le:relprime} we also deduce the following corollary.
\begin{corollary}%
\label{cor:h_K}%
\[%
\dd_{\calW_K}(x_{lm_1})  \Big{|}\nn_{\calW_K} \left (\frac{x^{h_K}_{lm_1}}{x^{h_K}_{klm_1}}-k^{2h_K}\right)^2
\]%
\end{corollary}%
\begin{proof}
From an elementary algebra calculation we have
\[
\frac{x^{h_K}_{lm_1}}{x^{h_K}_{klm_1}}-k^{2h_K}=\left(\frac{x_{lm_1}}{x_{klm_1}}-k^{2}\right)\sum_{r=0}^{h_K-1} \left(\frac{x_{lm_1}}{x_{klm_1}}\right)^{h_K-1-r}k^{2r},
\]
and therefore
\[
\nn_{\calW_K} \left(\frac{x_{lm_1}}{x_{klm_1}}-k^{2}\right)\Big{|}  \nn_{\calW_K} (\frac{x^{h_K}_{lm_1}}{x^{h_K}_{klm_1}}-k^{2h_K})\dd_{\calW_K} (\sum_{r=0}^{h_K-1} \left(\frac{x_{lm_1}}{x_{klm_1}}\right)^{h_K-1-r}k^{2r}).
\]
However, the only primes which can appear in
\[
\dd_{\calW_K} (\sum_{r=0}^{h_K-1} \left(\frac{x_{lm_1}}{x_{klm_1}}\right)^{h_K-1-r}k^{2r})
\]
are the primes occurring in
\[
\dd_{\calW_K} \left(\frac{x_{lm_1}}{x_{klm_1}}\right).
\]
The non-inverted part of the  divisor of $\displaystyle \frac{x_{lm_1}}{x_{klm_1}}$ is equal to  $\displaystyle \frac{\nn_{\calW_K} (x_{lm_1})\dd_{\calW_K} (x_{klm_1})}{\nn_{\calW_K} (x_{klm_1})\dd_{\calW_K} (x_{lm_1})}$, where $\displaystyle\frac{\dd_{\calW_K} (x_{klm_1})}{\dd_{\calW_K} (x_{lm_1})}$ is an integral divisor by  Lemma \ref{le:orderchangeold}.  This leaves only primes from $\nn_{\calW_K} (x_{klm_1})$ in the denominator.  Since none of these primes is present in  $\dd_{\calW_K} (x_{lm_1})$ due to Lemma \ref{le:relprime}, we have that
\[
\begin{array}{c}
\dd_{\calW_K} (x_{lm_1}) \Big{|} \nn_{\calW_K} \left (\frac{x^{h_K}_{lm_1}}{x^{h_K}_{klm_1}}-k^{2h_K}\right)^2\\
\Big \Updownarrow\\
\dd_{\calW_K} (x_{lm_1}) \Big{|} \nn_{\calW_K}\left [ \left(\frac{x_{lm_1}}{x_{klm_1}}-k^{2}\right)^2\left(\sum_{r=0}^{h_K-1} \left(\frac{x_{lm_1}}{x_{klm_1}}\right)^{h_K-1-r}k^{2r}\right)^2\right ]\\
\Big \Updownarrow\\
 \dd_{\calW_K} (x_{lm_1}) \Big{|} \nn_{\calW_K}\left(\frac{x_{lm_1}}{x_{klm_1}}-k^2\right)^2
 \end{array}
\]
\end{proof}
\begin{lemma}%
\label{le:square}
For any $k \in \Z_{>0}$ we have that $\dd(x_k)$, $\dd_k$ are squares of some integral divisors of $K$.
\end{lemma}
\begin{proof}%
From the Weierstrass equation $y^2=x^3+ax+b$ we have that for any prime $\pp$ of $K$, if $\ord_{\pp}x <0$, then $\ord_{\pp}(x^3+ax+b) = \ord_{\pp}x^3 <0$ and  $\ord_{\pp}y <0$ implying that $\ord_{\pp}x \equiv 0 \mod 2$.
\end{proof}%

\section{Diophantine Definition of Multiplication on Indices}
\label{sec:diophdef}
\setcounter{equation}{0}
We start with a basic fact and some easy lemmas.
\begin{lemma}
\label{le:notzero}
The set $\{x \in O_{K,\calW_K}: x \not=0\}$ is Diophantine over $O_{K,\calW_K}$.  (See Definition 2.2.3 and Proposition 2.2.4 of \cite{Sh34}.)
\end{lemma}
We now use the fact that we can define the set of non-zero integers of our ring to define relative primality over the ring.
\begin{lemma}%
\label{le:reprime}%
The set ${\tt R}=\{(A,B) \in O_{K,\calW_K}^2: AB \not = 0 \land (A,B)_{\calW_K}=1 \}$ is Diophantine over  $O_{K,\calW_K}$.
\end{lemma}%
\begin{proof}%
It is easy to see with the help of the Strong Approximation Theorem that  for $$(A,B) \in O_{K,\calW_K}^2$$ with $AB \not = 0$ the following statements are equivalent
\be  %
\item $(A,B)_{\calW_K}=1$
\item $\exists X, Y \in O_{K,\calW_K}: XA + YB = 1$
\ee%
\end{proof}%
\begin{notation}%
\label{not:2}
We define three sets: one to represent the points on our elliptic curve,  one to represent the elliptic curve addition, and one to represent the divisors of the denominators:
\begin{enumerate}
\item Let
\[%
{\tt E}= \{(U, V, X, Y) \in O_{K,\calW_K}^4 \mid  \exists k \in \Z_{\not = 0}:\frac{U}{V}=x_{m_0k}, \frac{X}{Y}=y_{m_0k}\}.
\]%
For each quadruple $(U, V, X, Y)$ the index $k=k(U,V,X,Y)$ will be unique (since the size of the torsion group divides $m_0$) and will be called the corresponding (to $(U,V,X,Y)$) index.
\item Let
\[%
{\tt Plus} = \{(U_1,V_1,X_1,Y_1), (U_2, V_2,X_2,Y_2), (U_3, V_3,X_3,Y_3)\} \subset {\tt E}^3
\]%
consist of triples of quadruples possessing corresponding indices $k_1, k_2, k_3$ satisfying $$k_1+k_2=k_3.$$

\item \label{it:d} Given $(U,V,X,Y) \in \tt E$, let
\[
d(U,V,X,Y) =\{(A,B) \in O_{K,\calW_K}^2:  \left (\frac{U}{V}\right )^{h_K} =\frac{A}{B}, (A,B)_{\calW_K} = 1\}.
\]
\end{enumerate}
\end{notation}

\begin{remark}
\label{rem:d}
The reason for defining the set $d(U,V,X,Y)$ is that over an arbitrary number field $K$ we cannot make sure that the numerators and denominators are relatively prime in our ring.  Thus a denominator can have ``too many'' primes in it and the divisibility conditions from Proposition \ref{prop:canmultiply} can fail if we replace the divisors by the denominators.  At the same time, by the definition of the class number, if we raise the $x$-coordinate to the power equal to the class number, we can obtain a relatively prime numerator and denominator.
\end{remark}

Given Lemma \ref{le:reprime}, the following assertion is obvious.
\begin{lemma}%
${\tt E}$, ${\tt Plus}$, and $d(U,V,X,Y)$ for  fixed values of $U, V, X, Y$, are Diophantine over $O_{K,\calW_K}$.
\end{lemma}

The next lemma and its corollary establish a connection between $ d(U,V,X,Y) $ and the divisor $\ff_k$ of the corresponding point on the elliptic curve.
\begin{lemma}%
\label{le:denom}
If $(U, V, X, Y) \in \tt E$,  $(A,B) \in d(U,V,X,Y)$, and $k$ is the corresponding index, then for all $\pp \not \in \calW_K$ we have that $h_K\ord_{\pp}\ff_k=\ord_{\pp}\nn_{\calW_K}(B)$ (Here we remind the reader that $\nn_{\calW_K}(B)$ is the non-inverted part of the numerator of the divisor of $B$).
\end{lemma}%
\begin{proof}%
By definition of ${\tt E}$ and  $d(U,V,X,Y)$ we have that  $\displaystyle \frac{A}{B} = x_{m_0k}^{h_K}$ for the corresponding to $(U,V,X,Y)$ index $k \in \Z_{\not=0}$.  Without loss of generality we can assume that $k >0$.  (``$-k$'' gives the same $B$ and the same $\ff_k$ by Lemma \ref{le:minus}.)  Let $\pp \not \in \calW_K$ be such that $\ord_{\pp}x_{m_0k} <0$.  Then either $\ord_{\pp}A <0$ or $\ord_{\pp}B >0$.  The first alternative is impossible because $A \in O_{K,\calW_K}$ and $\pp \not \in \calW_K$.  Hence we conclude that $\ord_{\pp}B >0$.  Further we also have that $\ord_{\pp}A=0$ because otherwise the relative primeness conditions requiring that $A$ and $B$ are not simultaneously divisible by any prime outside $\calW_K$ are violated.  Now we see that $$h_K\ord_{\pp}x_{m_0k} = \ord_{\pp}A-\ord_{\pp}B=-\ord_{\pp}B.$$  Suppose now that for some $\pp \not \in \calW_K$ it is the case that  $\ord_{\pp}x_{m_0k} \geq 0$ and $\ord_{\pp}B >0$.  In this case we also must have that $\ord_{\pp}A >0$ which again is impossible since $(A,B)_{\calW_K}=1$.

\end{proof}%

Given the lemma above we immediately conclude the following.
\begin{corollary}%
\label{cor:BC}
 If $I, I_1 \subset I, I_2 \subset I$ are finite subsets of non-zero integers,  $(U_i,V_i,X_i,Y_i) \in {\tt E}, (A_i,B_i) \in d(U_i,V_i,X_i,Y_i),  i \in I$ with $k_i$ being the corresponding indices,  then
  \[%
  \left (\prod_{i \in I_1}B_i\right ) \Big{|}_{\calW_K} \left (\prod_{i \in I_2 }B_i\right ) \iff \left(\prod_{i \in I_1}\ff_{k_i} \right )\Big{|} \left (\prod_{i \in I_2} \ff_{k_i}\right ).
  \]%
\end{corollary}%
Next we show that divisibility of indices is Diophantine in our ring.
\begin{lemma}
\label{le:divisibility}
If ${\tt Divide} = \{(U_1,V_1,X_1,Y_1), (U_2, V_2,X_2,Y_2)\} \subset {\tt E}^2$ consists of pairs of quadruples with the corresponding indices $k_1$ and $k_2$ such that $k_1 | k_2$, then {\tt Divide} is Diophantine over $O_{K,\calW_K}$.
\end{lemma}
\begin{proof}
If $(A_i,B_i) \in d(U_i,V_i,X_i,Y_i), i=1,2$, then by Corollary \ref{cor:BC} we have that $B_1 \Big{|}_{\calW_K} B_2$ if and only if $\ff_{k_1} | \ff_{k_2}$.  At the same time by Corollary \ref{cor:div2}, Part 3 and Lemma \ref{le:minus}, we have that  $\ff_{k_1} | \ff_{k_2}$ if and only if $k_1 | k_2$.
\end{proof}
We can now define multiplication on the {\it absolute values} of indices.

\begin{lemma}%
\label{le:newproduct}
Let $(U_i,V_i,X_i,Y_i) \in {\tt E}, (A_i, B_i) \in d(X_i,Y_i,U_,V_i), i=1,2,3$ with
\begin{equation}
\label{eq:relprime}
(B_1,B_2)_{\calW_K}=1
\end{equation}
 and
\begin{equation}%
\label{eq:product}
B_1B_2 \Big{|}_{\calW_K} B_3 \mbox{ and } B_3 \Big{|}_{\calW_K}  B_1B_2.
\end{equation}%
Then for the corresponding indices $k_1, k_2, k_3 \in \Z_{\not = 0}$ we have that $|k_1||k_2|=|k_3|$.
\end{lemma}%
\begin{proof}%
If $k_i$ is the index corresponding to   $(U_i,V_i,X_i,Y_i)$, then from (\ref{eq:relprime}) and  Lemma \ref{le:denom} we conclude that $(\ff_{k_1},\ff_{k_2})=1$.  Now from  Corollary
\ref{cor:BC} and (\ref{eq:product}) it follows that $\ff_{k_1}\ff_{k_2} = \ff_{k_3}$, and the assertion of the lemma is true by Lemma \ref{le:product}.
\end{proof}%

Our final step in this section is to define a square of an index.  This is all we need to define multiplication.
\begin{lemma}%
\label{le:define}
 Let  $(U_1, V_1,X_1,Y_1) \in {\tt E}$ be given and let $k_1$ be the corresponding index.  Assume there  exist quadruples $(U_i,V_i,X_i,Y_i) \in {\tt E}$ with the corresponding indices $k_i$ for $i=2, \ldots,8$  such that the following conditions and equations are satisfied.

\begin{equation}
\label{eq:B}
(A_i, B_i) \in d(U_i,V_i,X_i,Y_i), i=1,\ldots,8
\end{equation}
\begin{equation}%
\label{eq:mod4}%
 k_1 \equiv 4 \mod 16
\end{equation}
\begin{equation}%
\label{eq:mod16}%
 k_2 \equiv 1 \mod 2
\end{equation}
\begin{equation}%
\label{eq:1}
 k_3=k_1+ k_2
\end{equation}%
\begin{equation}%
\label{eq:2}
(B_1, B_2)_{\calW_K}= (B_1, B_3)_{\calW_K}=1
\end{equation}%
\begin{equation}%
\label{eq:3}
B_1B_2 \Big{|}_{\calW_K} B_4 \mbox{  and }   B_4 \Big{|}_{\calW_K}B_1 B_2,
\end{equation}%
\begin{equation}%
\label{eq:4}
B_1B_3 \Big{|}_{\calW_K} B_5 \mbox{ and }   B_5 \Big{|}_{\calW_K}B_1 B_3
\end{equation}%
\begin{equation}%
\label{eq:5}
 k_6=k_5-k_4
\end{equation}
\begin{equation}%
\label{eq:6}
 k_6 \equiv 0 \mod 16
\end{equation}
\begin{equation}%
\label{eq:7}
 k_7=k_1-1
\end{equation}
\begin{equation}%
\label{eq:8}
 k_8=k_6-1
\end{equation}
\begin{equation}%
\label{eq:9}
B_7 \Big{|}_{\calW_K} B_8,
\end{equation}%
Then $k_6=k_1^2$.  Conversely, if $k_1 \equiv 4 \mod 16$ then there exist  $(U_i,V_i,X_i,Y_i) \in {\tt E}, i=2, \ldots,8$
such that all the equations and conditions above can be satisfied.
\end{lemma}%

\begin{proof}%
First assume that for some $(U_i,V_i,X_i,Z_i) \in {\tt E}, i=1, \ldots,8$ with the corresponding indices $k_1,\ldots, k_8$ respectively, the conditions and equations (\ref{eq:B}) -- (\ref{eq:9}) are satisfied. Then by Lemma \ref{le:newproduct} we have from equations (\ref{eq:2}) -- (\ref{eq:4}) that $|k_4| =|k_1k_2|$ and $|k_5|=|k_1(k_1+k_2)|$. Thus,
$k_6 =\pm k_1^2$ or $k_6=\pm(k_1^2+2k_1k_2)$. From equations (\ref{eq:mod4}) and (\ref{eq:mod16}) we know that $k_1 \equiv 4 \mod
16$ and $k_2$ is odd.  Therefore, $k^2_1+2k_1k_2\not \equiv 0 \mod 16$. Thus, we must conclude that $k_6 =\pm k_1^2$.  Finally,
if $k_6 = -k_1^2$, then $k_8 = -1-k_1^2$ and consequently $k_1-1$ does not divide $k_8$, since $|k_1|\geq 4$, implying by Corollary \ref{cor:div2}
that $\ff_{(k_1-1)}$ does not divide   $\ff_{k_8}$.  Thus if $k_6=-k_1^2$, then (\ref{eq:9}) cannot hold.

Assume now that     $k_1 \equiv 4 \mod 16$.    By Proposition \ref{prop:directly} we can find a  $w \in
\Z_{\mbox{odd}}$ so that pairs $(k_1,k_2=w)$  and  $(k_1,k_3=k_1+k_2)$ can be multiplied
directly.     Let $k_4=k_1k_2, k_5 = k_1k_3=k_1^2 + k_1k_2$.  Let $k_6=k_5-k_4=k_1^2, k_7=k_1-1$,
and finally $k_8=k^2-1$ and define $(U_i,V_i,X_i,Y_i), i=2,\ldots,8$ using the definition of $\tt E$,  and $B_i,i=1,\ldots,8$ using the definition of $d(U_i,V_i,X_i,Y_i)$.  This will satisfy \eqref{eq:B}.  Note further that  Equations (\ref{eq:mod4}) and (\ref{eq:mod16}) can be satisfied by the choice of $k_1$ and $k_2$.
 Equation (\ref{eq:2}) will be satisfied by the definition of     ``can be multiplied directly''.
Equations (\ref{eq:3}) and (\ref{eq:4}) will be satisfied by  Proposition \ref{prop:canmultiply} and by the definition of     ``can be multiplied directly''.
  Equations (\ref{eq:5}) -- (\ref{eq:8}) will be satisfied by construction.  Finally Equation
(\ref{eq:9}) will be satisfied by Corollary \ref{cor:div2}.
\end{proof}%

Lemma \ref{le:define} completes the proof of Theorem \ref{thm:main} and Corollary \ref{cor:main}.  (The density computation is in the Appendix.)\\

We finish this section with a new notation to be used below.

\begin{notation}%
\label{not:pi}
\begin{itemize}
\item Given $(U_i,V_i,X_i,Y_i) \in {\tt E}, i =1,2,3$ we will say that $$((U_1,V_1,X_1,Y_1), (U_2,Y_2,X_2,Y_2), (U_3,Y_3,X_3,Y_3) ) \in \Pi$$ to mean that the corresponding indices $k_1, k_2, k_3$ satisfy $k_3=k_1k_2$.
\item Let
\[%
{\tt E}_1= \{(U, V, X, Y) \in O_{K,\calW_K}^4 \mid  \exists \mbox{ unique } k \in \Z_{\not = 0}:\frac{U}{V}=x_{m_1m_0k}, \frac{X}{Y}=y_{m_1m_0k}\}.
\]%
The positive integer $m_1$ is defined in Lemma \ref{le:equiv}.
\end{itemize}
\end{notation}

\section{Defining $\Z$ over $O_{K,\calW_K}$ Using One Universal Quantifier}
\setcounter{equation}{0}
In this section we use the existential definition of multiplication on indices to give a first-order definition of $\Z$ over $O_{K,\calW_K}$ using just one universal quantifier.  We start with a technical lemma.
\begin{lemma}%
\label{le:firstorder}
 If $z \in O_{K,\calW_K}$ has the following property:   $$ \exists  U_1, V_1, X_1, Y_1, \, \forall b, \, \exists U_2, V_2, X_2, Y_2, U_3, V_3, X_3, Y_3, A_1, A_2, A_3, B_1, B_2, B_3,C$$ (with all the variables ranging over $O_{K,\calW_K}$) such that

\begin{equation}
\label{eq:multiply}
(U_1,V_1,X_1,Y_1), (U_3,V_3,X_3,Y_3) \in {\tt E}, (U_2,V_2,X_2,Y_2) \in {\tt E}_1,
\end{equation}
\begin{equation}
\label{eq:multiply2}
((U_1,V_1,X_1,Y_1), (U_2,V_2,X_2,Y_2), (U_3,V_3,X_3,Y_3)) \in \Pi,
\end{equation}

\begin{equation}
\label{eq:E}
(A_i, B_i) \in d(U_i,V_i,X_i,Y_i), i=1,2,3,
\end{equation}
\begin{equation}
\label{eq:divide}
b^{2h_K} \mid_{\calW_K} B_2,
\end{equation}%

\begin{equation}%
\label{eq:congruence}
(A_3B_2z - B_3A_2)^{2h_K} = B_2^{2h_K+1}C,
 \end{equation}
then $z \in \Z$.

Conversely, if $z_0 \in \Z_{\not = 0}$ and $z=z_0^{2h_K}$, then Equations \eqref{eq:multiply}--\eqref{eq:congruence} can be satisfied with variables as above ranging over $O_{K,\calW_K}$.
 \end{lemma}%
\begin{proof}%
First assume that equations above are satisfied for some $z \in O_{K,\calW_K}$. From  (\ref{eq:multiply})--\eqref{eq:E},  we conclude that if $k_1, k_2, k_3$ are the indices corresponding to 
\[
(U_1, V_1, X_1, Y_1), (U_2, V_2, X_2, Y_2), \mbox{ and }(U_3, V_3, X_3, Y_3)
\]
 respectively, then $k_3=k_1k_2$,  $k_2 \equiv 0 \mod m_1$ and  $\nn_{\calW_K}(B_i) =\ff^{h_K}_{k_i}$.  Further for the discussion below $k_1$ is fixed.  From equation (\ref{eq:congruence}), we obtain that
\[%
\nn_{\calW_K}(B_2^{2h_{K}+1}) \Big{|} \nn_{\calW_K}(A_3B_2z -A_2B_3)^{2h_K}
\]%
and therefore
\[%
\nn_{\calW_K}(B_2) \Big{|} \nn_{\calW_K}(A_3z -\frac{A_2B_3}{B_2})^{2h_K}
\]%
Further, since $(B_2,A_3)_{\calW_K}=1$ by Lemma \ref{le:relprime}, we have that
\begin{equation}%
\label{eq:need}
\nn_{\calW_K}(B_2) \Big{|} \nn_{\calW_K}(z -\frac{A_2B_3}{A_3B_2})^{2h_K}.
\end{equation}%
Thus,  since $\nn_{\calW_K}(B_2)$ is a $2h_K$-th power of another divisor in $K$ by Lemma \ref{le:square},  and by the definition of $B_2$ we have that
\begin{equation}
\label{eq:start}
\sqrt[2h_K]{\nn_{\calW_K}(B_2)}  \Big{|}_{\calW_K} \nn_{\calW_K}(z -\frac{A_2B_3}{A_3B_2}).
\end{equation}
From Corollary \ref{cor:h_K} and Lemma \ref{le:square}, since $k_2$ is divisible by $m_1$ while $k_3=k_1k_2$, and therefore $k_3$ is also divisible by $m_1$, we conclude that
\[
\sqrt{\dd_{\calW_K}(x_{k_2m_0})} \Big{|}_{\calW_K} \nn_{\calW_K}\left (\frac{x^{h_K}_{k_2m_0}}{x^{h_K}_{k_3m_0}}-k_1^{h_K}\right ),
\]
and therefore, using the definition of $B_2$, we have
\[
 \sqrt[2h_K]{\nn_{\calW_K}(B_2)}\Big{|}_{\calW_K}  \nn_{\calW_K}\left(\frac{x^{h_K}_{k_2m_0}}{x^{h_K}_{k_3m_0}}-k_1^{h_K}\right ).
\]
Substituting $\displaystyle \frac{A_2B_3}{A_3B_2}$ for $\displaystyle \frac{x^{h_K}_{k_2m_0}}{x^{h_K}_{k_3m_0}}$ we obtain
\begin{equation}
\label{eq:subst}
 \sqrt[2h_K]{\nn_{\calW_K}(B_2)}\Big{|}_{\calW_K}  \nn_{\calW_K}\left (\frac{A_2B_3}{A_3B_2}-k_1^{h_K}\right ).
\end{equation}
Combining \eqref{eq:start} and \eqref{eq:subst} and using the definition of the divisibility in the divisor semigroup, we obtain
\[%
\sqrt[2h_K]{\nn_{\calW_K}(B_2)}  \Big{|} \nn_{\calW_K}(z -k_1^{2h_K}),
\]%
and
\[%
\nn_{\calW_K}(b)  \Big{|} \nn_{\calW_K}(z -k_1^{2h_K}),
\]%

Since the last divisibility condition has to hold for all $b$,  we must conclude that $z=k_1^{2h_K}$.\\

Conversely, suppose $z=z_0^{2h_K}$ for $z_0 \in \Z_{\not=0}$. Let $(U_1,V_1, X_1,Y_1) \in \tt E$ with the corresponding index $k_1=z_0$.  Let $b \in O_{K,\calW_K}$ be given.  Let $k_2 \equiv 0 \mod m_1$ be such that $b^2 \mid_{\calW_K} \dd_{\calW_K}(x_{k_2m_0})$.  Such  an index $k_2$ exists by Lemma \ref{le:anydivisor}.  Let   $(U_2,V_2, X_2,Y_2) \in \tt E$ correspond to $k_2$.  Let $k_3=k_1k_2$  and let $(U_3, V_3, X_3,Z_3) \in \tt E$ correspond to the index $k_3$.   Observe that conditions (\ref{eq:multiply}) and (\ref{eq:multiply2}) are now satisfied. Further note that equation (\ref{eq:need}) holds by Corollary \ref{cor:h_K} and therefore
equation (\ref{eq:congruence}) holds also.
\end{proof}%
To deal with the case of an arbitrary non-zero integer we add the following corollary.
\begin{corollary}
\label{cor:firstorder}%
 If $z_0 \in O_{K,\calW_K}$ has the following property: $$\exists z_1,\ldots, z_{2h_K},$$ $$ \exists  U_{1,0}, \ldots, U_{1,2h_K}, V_{1,0}, \ldots, V_{1,2h_K},$$ $$ \exists X_{1,0}, \ldots,  X_{1,2h_K}, Y_{1,0}, \ldots,  Y_{1,2h_K},$$ $$ \forall b,$$  $$\exists U_2, V_2, X_2, Y_2,$$ $$\exists U_{3,0}\ldots  U_{3,2h_K}, V_{3,0}, \ldots, V_{3,2h_K},$$  $$\exists X_{3,0}, \ldots,  X_{3,2h_K}, Y_{3,0}, \ldots,  Y_{3,2h_K},$$  $$\exists A_{1,0}, B_{1,0}, \ldots, A_{1,2h_K}, B_{1,2h_K}, A_2, B_2, A_{3,0}, B_{3,0}, \ldots, A_{3,2h_K}, B_{3,2h_K}, C_0,\ldots,C_{2h_K}, $$ (with all the variables ranging over $O_{K,\calW_K}$) such that
 \begin{equation}
 \label{eq:z}
 z_j = (z_0 +j)^{2h_K}, j=0,\ldots, 2h_K,
 \end{equation}
\begin{equation}
\label{eq:multiply3}
(U_{i,j},V_{i,j},X_{i,j},Y_{i,j}) \in {\tt E} , i =1,3, j=0,\ldots, 2h_K,
\end{equation}
\begin{equation}
\label{eq:multiply11}
(U_2,V_2,X_2,Y_2) \in {\tt E}_1,
\end{equation}
\begin{equation}
\label{eq:multiply31}
[(U_{1,j},V_{1,j},X_{1,j},Y_{1,j}), (U_2,V_2,X_2,Y_2), (U_{3,j},V_{3,j},X_{3,j},Y_{3,j})] \in \Pi, j=0,\ldots, 2h_K,
\end{equation}
\begin{equation}
\label{eq:E1}
(A_{i,j},B_{i,j}) \in d(U_{i,j},V_{i,j},X_{i,j},Y_{i,j}), i=1,3, j=0,\ldots,2h_K,
\end{equation}
\begin{equation}
\label{eq:E2}
(A_2,B_2) \in d(U_2,V_2,X_2,Y_2),
\end{equation}
\begin{equation}
\label{eq:divide1}
b^{2h_K} \mid_{\calW_K} B_2,
\end{equation}%
\begin{equation}%
\label{eq:congruence2}
(A_{3,j}B_2z_j - B_{3,j}A_2)^{2h_K} = B_2^{2h_K+1}C_j, j=0,\ldots,2h_K,
 \end{equation}
then $z_0 \in \Z$.

Conversely, if $z_0 \in \Z_{\not = 0}$, then Equations \eqref{eq:z} -- \eqref{eq:congruence2} can be satisfied with variables as described above ranging over $O_{K,\calW_K}$.

\end{corollary}%
\begin{proof}%
If the assumptions of the corollary are true, then by Lemma \ref{le:firstorder} we have that
$$z_0^{2h_K},\ldots,(z_0+2h_K)^{2h_K} \in \Z,$$ and by Corollary B.10.10 of \cite{Sh34}, we have that $z_0 \in \Q$.  At the
same time, since $z_0^{2h_K} \in \Z$ we have that   $z_0$ is an algebraic  integer, and hence in $\Z$.
The rest of the proof is analogous to the proof of the second part  of Lemma \ref{le:firstorder}
\end{proof}%

The last proposition concludes the proof of Theorem \ref{thm:first-order}.
\section{Infinite Extensions}
\begin{notationassumptions}%
We add the following to our assumption list.
\begin{itemize}%
\item Let $K_{\infty}$ be a possibly infinite algebraic extension of $K$.
\item Assume $E(K_{\infty})=E(K)$.
\item Let $O_{K_{\infty},\calW_{K_{\infty}}}$ be the integral closure of $O_{K,\calW_K}$ in $K_{\infty}$.
\end{itemize}%
\end{notationassumptions}%
Given the assumptions on our elliptic curve, it is easy to see that the results of the previous section will carry over, and
therefore we have the following theorem:
\begin{theorem}%
\label{th:infinite}
\be%
\item Let $K$ be a number field.  Let $E$ be an elliptic curve defined and of rank one over $K$.  Let $P$ be a
generator of $E(K)$ modulo the torsion subgroup, and fix an affine Weierstrass equation for $E$ of the form $y^2=x^3+ ax +b$, with $a, b \in O_K$, where $O_K$ is the ring of integers of $K$.  Let $(x_n,y_n)$ be the coordinates of $[n]P$ derived from this Weierstrass equation.  Then there exists a set of $K$-primes $\calW_K$ of natural density one,  and a positive integer $m_0$ such that the following set $\Pi_{\infty} \subset O_{K_{\infty},\calW_{K_{\infty}}}^{12}$ is Diophantine over $O_{K_{\infty},\calW_{K_{\infty}}}$.
\[%
(U_1, U_2, U_3, X_1,X_2,X_3, V_1,V_2, V_3, Y_1, Y_2, Y_3)\in \Pi_{\infty} \Leftrightarrow
\]%
\[%
\exists \mbox{ unique } k_1, k_2, k_3 \in \Z_{\not= 0}   \mbox{ such that } \left(\frac{U_i}{V_i}, \frac{X_i}{Y_i}\right )=(x_{m_0k_i}, y_{m_0k_i}) \mbox{ and } k_3=k_1k_2.
\]%

\item For $n \not = 0$ let $\phi_{\infty}(n)=[(U_n,V_n, X_n,Y_n)]$, the class of  $(U_n,V_n, X_{n}, Y_{n})$ under the equivalence relation described below, where $U_n, V_n, X_n,Y_n \in
O_{K_{\infty},\calW_{K_{\infty}}}$, $Y_nV_n \not = 0$, and $\displaystyle (x_{m_0n},y_{m_0n})=\left (\frac{U_n}{V_n}, \frac{X_n}{Y_n}\right)$.  Let $\phi_{\infty}(0)=\{[0,0,0,0]\}$.  Then
$\phi_{\infty}$  is a class Diophantine model of $\Z$.  (Here if $YV \not =0$ we have that  $(U,V, X,Y) \approx (\hat U,\hat V, \hat X,\hat Y)$ if and only if  $\displaystyle (\frac{\hat U}{\hat V}, \frac{\hat X}{\hat Y})=(\frac{U}{V}, \frac{X}{Y}).$)

\item $\Z$ is definable over $O_{K_{\infty},\calW_{K_{\infty}}}$ using one universal quantifier.
\end{enumerate}
\end{theorem}%

\section{Appendix}
\label{sec:density}
\setcounter{equation}{0}
In this Appendix we calculate the natural density of $\calV_K$.  This calculation is similar to the one carried in \cite{CS}.  We use Notation \ref{S:notation section} and a new notation: for a prime $\pp$ of a number field $K$ we let ${\mathbf N\pp}$ denote the size of the residue field of $\pp$.

\begin{lemma}%
\label{le:pindenom}%
Let $\ell \in \calP(\Q)$ and suppose $\pp \in \calS_{\ell^{n+1}}\setminus {\calS}_{\ell^n}$ for some $n \in \Z_{\geq 0}$.   (Such a $\pp$ exists, if
$n\geq a_\ell$.) Then $\ell^{n+1} < 3\norm \pp$.
\end{lemma}%
\begin{proof}%
 If $\pp \in {\calS}_{\ell^{n+1}}\setminus {\calS}_{\ell^n}$, then  $\pp$ does not divide the discriminant of our Weierstrass equation and $\tilde{E}$, the reduction of $E \mod \pp$  is non-singular.
Further, $x_{\ell^n}$, $y_{\ell^n}$ are integral at $\pp$, while $\ord_{\pp}x_{\ell^{n+1}}<0$, $\ord_{\pp}y_{\ell^{n+1}}<0$. Therefore, under the reduction mod  $\pp$, the image of $[\ell^n]P$ is not $\tilde{O}$ -- the image of $O$ mod $\pp$ , while $[\ell^{n+1}]\tilde{P}=\tilde{O}$. Thus we must conclude that $E(\F_{\pp})$ has
an element of order $\ell^{n+1}$ and therefore $\ell^{n+1}| \#E(\F_{\pp})$.  Let $\# \F_{\pp}= \norm{\pp} = q$. From
a theorem of Hasse we know that $\#E(\F_{\pp}) \leq q+1+2 \sqrt{q} \leq 3q$ (see \cite{Sil1}, Chapter V,
Section 1, Theorem 1.1).
\end{proof}%
    \begin{lemma}%
\label{prop:density}
The natural density of the set  $\calA=\{\pp_{\ell^k}: \ell \in \calP_{\Q}, k \in \Z_{>1} \land k \geq
a_{\ell}\}$ is zero.
\end{lemma}%
\begin{proof}%
For $\pp = \pp_{\ell^k} \in \calA$, the preceding lemma says that $3\norm\pp_{\ell^k} > \ell^k$.  Thus, since each $\pp \in \calA$ corresponds to a distinct pair $(\ell, k)$  with $\ell \in \calP(\Q)$ and $k \in \Z_{\geq 2}$ with $3\norm\pp > \ell^k$, we have the following inequality:
\[%
\# \{\pp \in  \calA \, : \, \norm{\pp} \leq X \}  \leq   \# \{(\ell,k)  \in \calP_\Q \times \Z_{k\geq 2} \, : \,  \ell \leq
\sqrt[k]{3X}\}
\]%
Clearly if $\sqrt[k]{3X} < 2$, there will be no prime $\ell$ with $\ell \leq \sqrt[k]{3X}$.  Thus,
we can limit ourselves to positive integers $k$ such that $k \leq \log_2 (3X)$.

By the Prime Number Theorem (see \cite{L}, Theorem 4, Section 5, Chapter XV), for some positive
constant $C$ we have that $\#\{\ell \in \calP_{\Q}: \ell \leq X\} \leq C {X}/{\log X} $ for all
$X \in \Z_{>0}$.  From the discussion above we now have the following sequence of inequalities:
\begin{eqnarray*}
\{\pp \in  \calA \, : \, \norm{\pp} \leq X \} &\leq& \sum_{k=2}^{\lceil \log_2 (3X) \rceil} \#\{\ell \in \calP_{\Q}: \ell \leq \sqrt[k]{3X}\}  \\ &\leq&
\sum_{k=2}^{\lceil \log_2(3X) \rceil} \#\{\ell \in \calP_{\Q}: \ell \leq \sqrt{3X}\} \\ &\leq& \log_2(3X)[C\frac{\sqrt{3X}}{\log{\sqrt{3X}}}] = \tilde C\sqrt{X}
\end{eqnarray*}
for some positive constant $\tilde C$.  At the same time by  the Prime Number Theorem again we
also know that  for some positive constant $\bar C$ we have
$
\#\{\pp \in \calP_K: \norm \pp \leq X \} \geq \bar C {X/\log X}.
$
Thus the upper density of $\calA$ is
\[%
\limsup_{X \rightarrow \infty}\frac{\#\{ \pp \in \calA \, : \, \norm \pp \leq X\}}{\#\{\pp \in
\calP_K: \norm \pp \leq X \}} \leq \limsup_{X \rightarrow \infty} \frac{\tilde C\sqrt{X}\log X}{\bar CX} = 0.
\]
Hence $\calA$ has a natural density, and it is zero.
\end{proof}

\begin{proposition}%
The set $\calV_K(P)$ has natural density zero.
\end{proposition}%
\begin{proof}%
We first observe  that it was proven in \cite{Po} and \cite{PS}  that the set
\[%
\calB=\{\pp_{\ell}: \ell \in \calP_{\Q} \land a_{\ell}=1\}
\]%
has a natural density that is zero. Finally we note that $\calB \cup \calA=\calV_K(P)$.
\end{proof}


\end{document}